\documentclass[preprint,groupedaddress,showkeys]{revtex4-1}
\usepackage{graphicx}
\usepackage{dcolumn}
\usepackage{bm}
\usepackage{amssymb}
\usepackage{amsthm}\usepackage{mathrsfs}
\usepackage{amsmath}
\usepackage{graphicx}
\usepackage{epstopdf}

\newtheorem{theorem}{\hskip\parindent\bf Theorem}
\newtheorem{lemma}{\hskip\parindent\bf Lemma}
\begin{document}

\preprint{}

\title[Bautin bifurcation in delayed reaction-diffusion systems]{Bautin bifurcation in delayed reaction-diffusion systems with application to the Segel-Jackson model}

\author{Yuxiao Guo, Ben Niu*}

\affiliation{Department of Mathematics, Harbin Institute of Technology, Weihai 264209,  China.\\*Corresponding author, niu@hit.edu.cn}

\date{\today}

\begin{abstract}
 In this paper, we present an algorithm for deriving the normal forms of Bautin bifurcations in reaction-diffusion systems with time delays and Neumann boundary conditions. On the center manifold near a Bautin bifurcation, the first and second Lyapunov coefficients are calculated explicitly, which completely determine the dynamical behavior near the bifurcation point. As an example, the Segel-Jackson predator-prey model is studied. Near the Bautin bifurcation we find the existence of  fold bifurcation of periodic orbits, as well as subcritical and supercritical Hopf bifurcations. Both theoretical and numerical results indicate that solutions with small (large) initial conditions converge to  stable periodic orbits (diverge to infinity).
\end{abstract}

                           \maketitle

\section {Introduction}

In recent years, Hopf bifurcation in reaction-diffusion systems with time delays has been widely studied \cite{Faria2001,Su-Wei-Shi}. The results are used to interpret the oscillation behaviors in many subjects, such as, biology \cite{Chen,Yi}, chemistry \cite{Wei}, and epidemics \cite{Jin}, etc. The authors usually  investigated the change of stability of equilibria and the existence  of Hopf bifurcations. Furthermore, they calculated the normal forms  near the Hopf bifurcations to determine the stability of bifurcating periodic solutions and direction of Hopf bifurcations. The framework of deriving normal forms is due to the center manifold reduction technique introduced in  \cite{Yi,Hassard,Wu,Lin}.

As is known to all, supercritical Hopf bifurcations with negative first Lyapunov coefficient lead a system to stable periodic orbits at the side of bifurcation point that the equilibrium becomes unstable; whereas subcritical Hopf bifurcations with positive first Lyapunov coefficient induce unstable periodic orbits at the side of bifurcation point that the equilibrium becomes stable. General theory about Hopf bifurcation can be found in \cite{Holmes,Wiggins2,Kuz} for ordinary differential equations, in \cite{Hassard} for functional differential equations, and in \cite{Wu,Henry} for partial differential equations or partial functional differential equations. If we introduce another parameter and with the change of that parameter, Hopf bifurcation points forms a curve in the two-parameter plane. Thus, at the intersection of the supercritical Hopf bifurcation and  subcritical Hopf bifurcation curves, there will be a point with zero first Lyapunov coefficient that makes possible Bautin bifurcation happen which is of codimension two \cite{Holmes,Wiggins2,Kuz,jmaa}. Near a Bautin bifurcation, the dynamical behavior of a system is quite complicated: the bifurcation point separates branches of subcritical and supercritical Hopf bifurcations in the parameter space. For  parameter values near the singularity, the system has two limit cycles which collide and disappear via a fold bifurcation of periodic orbits.

So far, there are few papers on  codimension two bifurcations in delay equation, most of which focus in retarded functional differential equations\cite{Guo, Jiang}, and neutral functional differential equations \cite{Niu1, Niu2}. The study on codimension two bifurcations of  delayed reaction-diffusion equations just begins in  recent  few years, and most of which are about  Turing-Hopf bifurcation \cite{Song}, Hopf-zero bifurcation\cite{An} or Double Hopf bifurcation\cite{Du}. Near these bifurcations the spatially non-homogeneous steady-state solutions have been found as well as periodic solutions and quasi-periodic solutions. Bautin bifurcations have been found and analyzed in ordinary differential equations \cite{ode1,ode2,ode3,ode4} or in delay differential equations \cite{Xu1,Xu2,Ion}. However, to our best knowledge, there are not any research papers found in literature about Bautin bifurcations in  reaction-diffusion systems with time delays.

Motivated by such a consideration, in this paper we will give  universal derivations of the first and second Lyapunov coefficients for Bautin bifurcations   in  reaction-diffusion equations with time delays, then the normal form of Bautin bifurcations is given and has two different kinds of bifurcation diagrams. The results are then applied to a Segel-Jackson prey-predator model\cite{Wang,Segel}
\begin{equation}\label{S_J_model1}
 \left\{
 \begin{array}{l}
 \frac{\partial u(x,t)}{\partial t}=d\Delta u(x,t)+[1+ku(x,t)]u(x,t)-au(x,t)v(x,t)\\
  \frac{\partial v(x,t)}{\partial t}=\Delta v(x,t)+u(x,t-\tau)v(x,t-\tau)-v^2(x,t)
 \end{array}~~~~~x\in[0,l\pi], t>0.
 \right.
 \end{equation}
 In such a system with time delay and homogeneous boundary conditions, we find that Bautin bifurcation appears in $(k,\tau)$ plane. Both subcritical and supercritical Hopf bifurcations are found theoretically near the bifurcation point.

The rest part of the paper will be organized as follows. In section 2,
the fifth-order normal form near a pair of imaginary roots in a general   reaction-diffusion system of $N$ equations with time delays are derived, from which detailed and explicit formulae to determine the first and second Lyapunov coefficients for Bautin bifurcation are given.  In section 3, we apply the results to a  Segel-Jackson model, in which we find a Bautin bifurcation, and both subcritical and supercritical Hopf bifurcations are theoretically analyzed. To support the results some simulations are performed to illustrate the dynamical behaviors near the Bautin bifurcation. Finally a conclusion section completes this paper.

 \section{Bautin bifurcation in a reaction-diffusion system with time delays}

 We consider the following general form of a reaction-diffusion system with delays,
\begin{equation}\label{abs1}\frac{d U(t)}{d t}=D(\mu)\Delta U(t)+L(\mu)(U_t)+F(\mu,U_t).\end{equation}
 Here the state variable  $U$ maps $[-1,+\infty)$ to $X$, and the state space is $$X=\left\{(u_1,u_2,\ldots,u_N)\in\left[{H^{2}(0,l\pi)}\right]^N: \left.\left(\frac{d u_1}{ dx},\frac{d u_2}{ dx},\ldots,\frac{d u_N}{ dx}\right)\right|_{x=0,l\pi}=0\right\}.$$
 Notice that here we use homogeneous Neumann boundary conditions, and a one dimensional space interval $[0,l\pi]$.
 Using the general notation of delayed equations \cite{Wu}, the state $U_t$ in phase space $\mathscr{C}:=C([-1,0],X)$ is defined by $U_t(\theta)=U(t+\theta)$ for any $\theta\in[-1,0]$. The bifurcation parameter $\mu=(\mu_1,\mu_2)\in\mathbb R^2$.
 The diffusion matrix $D(\mu)=\mathrm{diag}~(d_{1}(\mu),d_{2}(\mu)),\ldots,d_{n}(\mu))$.~
$L:{\mathbb{R}}^{2}\times\mathscr{C}\rightarrow X$ is a bounded linear operator.
$F:{\mathbb{R}}^{2}\times\mathscr{C}\rightarrow X$ is a $C^{k}$ $(k\geq5)$ function of the form $(F_1, F_2,\ldots,F_N)^{\rm T}$ and satisfies $F(\mu,0)=0$ and $D_{\varphi}F(\mu,0)=0$, which means $U_t=0$ is a steady states of system (\ref{abs1}).

The linearized equation of system (\ref{abs1}) at $U_t=0$ is
 \begin{equation}\label{absl}
 \frac{d U(t)}{d t}=D(\mu)\Delta U(t)+L(\mu)U_t.
 \end{equation}
Since  we are about to calculate the fifth order normal form at a Bautin bifurcation,
we first assume that when $\mu\in V$, with $V$ representing some neighbourhood of the origin in $\mathbb R^2$, the characteristic equation of (\ref{absl})
 \begin{equation}\label{cha}
 d e t\left[\lambda {\rm{I}}+\frac{n^2}{l^2}D(\mu)-L(\mu)(\rm  e^{\lambda\cdot}I)\right]=0
     \end{equation}
has a pair of roots $\alpha(\mu)\pm\rm  i\omega(\mu)$ for some $n=n_0$ where $\alpha(\mu)$ and $\omega(\mu)$ are both continuously differentiable, and $\alpha(0)=0$, $\omega(0)=\omega_0>0$. We further assume that all the other roots of (\ref{cha}) have negative real parts for $\mu\in V$. Thus the center manifold will be locally attractive and the local dynamical behaviors of (\ref{absl}) will be governed by the normal form on the center manifold \cite{Lin}.

The solution operator of Eq.(\ref{abs1}) generates a $C_0$-semigroup with the infinitesimal generator $A_\mu$ defined by

 \begin{equation}\label{inf2} A_\mu \phi=\left\{\begin{array}{ll} \dot\phi(\theta),&\theta\in [-1,0),\\
  D(\mu)\Delta\phi(0)+L(\mu)\phi,&\theta=0.\end{array}\right.\end{equation}
 The domain of $A_\mu$ is chosen by
 $${\rm dom}(A_\mu):=\{\phi\in\mathscr{C}:\dot{\phi}\in\mathscr{C}, \phi(0)\in {\rm dom}(\Delta), \dot{\phi}(0)=  D(\mu)\Delta\phi(0)+L(\mu)\phi\}.$$
We rewrite Eq.(\ref{abs1}) as an abstract ordinary differential equation
 \begin{equation}\label{abso}\frac{d}{dt}{U_t}=A_\mu U_t+R(\mu,U_t)\end{equation}
with the nonlinearity
 $$R(\mu, U_t)(\theta)=\left\{\begin{array}{ll}
      0,&   \theta\in[-1,0), \\
      F(\mu,U_t),& \theta=0.
 \end{array}\right.$$

Write the Fourier basis of $X$ by $\{b_n e_1,b_n e_2,\ldots,b_ne_N\}_{n=0}^{+\infty}$ with $e_1=(1,0,\ldots,0)^{\rm T}$,   $e_2=(0,1,\ldots,0)^{\rm T},\ldots,$ $e_N=(0,0,\ldots,1)^{\rm T}$, and $b_n=\frac{\cos \frac{nx}{l}}{\|\cos\frac{nx}{l}\|_{L_2}}$.
 Here the $L_2$ norm is $$\|\cos\frac{nx}{l}\|_{L_2}=\left\{\begin{array}{l}\sqrt{\frac{1}{l\pi}},~~n=0,\\
 \sqrt{\frac{2}{l\pi}},~~n\neq0.\end{array}\right.$$
 Also denote by $\beta_n=(b_n,b_n,\ldots,b_n)^{\rm T}$, then
 for any $\phi=(\phi^{(1)},\phi^{(2)},\ldots,\phi^{(N)})^{\rm T}\in \mathscr{C}$, we can define $\phi_n=<\phi,\beta_n>=(<\phi^{(1)},b_n>,<\phi^{(2)},b_n>,\ldots,<\phi^{(N)},b_n>)^{\rm T}$.

 Now we have the action of $A_\mu$ on $\phi_n b_n$, denoted by $A_{\mu,n}$ as
 $$A_{\mu,n}(\phi_nb_n)=\left\{\begin{array}{l}\dot\phi(\theta)b_n,\theta\in[-1,0),\\
 \int_{-1}^0 d\eta_n(\mu,\theta)\phi_n(\theta)b_n,\theta=0,
 \end{array}
 \right.$$
 with the function of bounded variation $\eta_n$ satisfying
 $$-\frac{n^2}{l^2}D(\mu)\phi_n(0)+L_\mu(\phi_n)=\int_{-1}^0 d\eta_n(\mu,\theta)\phi_n(\theta).$$

 Consider the adjoint operator of $A_0$, $A^\ast$ given by
 $$A^\ast\varphi(s)=\left\{\begin{array}{ll}
 -\dot\varphi(s),& s\in(0,1],\\
 \sum\limits_{n=0}^\infty\int_{-1}^0d\eta_n^\textrm{T}(0,t)\varphi_n(-t)b_n,& s=0,\end{array}
 \right.$$
 for some  $\varphi\in\mathscr{C}^\ast=\mathscr{C}([0,1],X)$ and $\varphi_n=<\varphi,b_n>$.
According to the formal adjoint theory given in \cite{Wu,Hale}, we define the bilinear form by
 \begin{equation}\label{varphi0}
     (\varphi,\phi)=\sum_{k,j=0}^\infty(\varphi_k,\phi_j)_c\int_0^{l\pi} b_k b_j dx.
 \end{equation}
 In fact, $\int_0^{l\pi} b_k b_j d x=0$, $\forall k\neq j$. Thus, we only need to define $(\varphi_n,\phi_n)_c$ from \cite{Hassard,Wu}, as
 $$(\varphi_n,\phi_n)_c=\bar\varphi_n^\textrm{T}(0)\phi_n(0)-\int_{-1}^0\int_{\xi=0}^\theta\bar\varphi_n^\textrm{T}(\xi-\theta)d\eta_n(0,\theta)\phi_n(\xi)d\xi.$$

 Let $q(\theta)b_{n_0}$ and $q^\ast(s)b_{n_0}$ be the eigenfunctions of $A_0$ and $A^\ast$ corresponding to $\rm  i\omega_0$ and $-\rm  i\omega_0$ respectively, with $({q^\ast},q)_c=1, (q^\ast,\bar{q})_c=0$.

 Write the centerspace
 $$P=\{z q b_{n_0}+\bar z\bar q b_{n_0}|z\in \mathbb C\}$$
and its orthogonal complement space
$$Q=\{\phi\in \mathscr{C}|(q^\ast b_{n_0},\phi)=0 ~and ~(\bar{q}^\ast b_{n_0},\phi)=0\},$$
then $\mathscr{C}=P\bigoplus Q$ and the state variable $U_t$ of Eq.(\ref{abso}) could be decomposed by
\begin{equation}\label{uuuu}U_t=z(t)q(\cdot)b_{n_0}+\bar{z}(t)\bar{q}(\cdot)b_{n_0}+W(t,\cdot)\end{equation}
where $W(t,\cdot)\in Q.$

Obviously
$z(t)=(q^\ast b_{n_0},U_t)$
and $W(t,\theta)=U_t(\theta)-2{\rm Re}\{z(t)q(\theta)b_{n_0}\}$.
Thus, we have
\begin{equation}\label{zdot22}\dot{z}(t)=\textrm{i}\omega_0 z(t)+\bar{q}^{\ast\textrm{T} }(0)<F(0,U_t),\beta_{n_0}>,
\end{equation}
where the nonlinear term is given by
$$<F,\beta_n>=(<F_1,b_n>,<F_2,b_n>,\ldots,<F_N,b_n>)^\textrm{T}.$$

Using the center manifold theorem \cite{Wu,Lin},  $W(t,\cdot)$ can be expressed by a series of $z$ and $\bar{z}$, starting at quadratic terms
\begin{equation}\label{cm22}W(t,\cdot)=W(z,\bar z,\cdot)=\sum\limits_{i+j\geq2}W_{ij}(\cdot)\frac{z^i\bar z^j}{i!j!}
\end{equation}
where $z$ and $\bar{z}$ are approximately the local coordinates for center manifold $\mathscr{C}_0$ in two directions $q b_{n_0}$ and $\bar{q}b_{n_0}$, respectively. For solution $U_t\in\mathscr{C}_0,$ we have $U_t=z(t)q(\cdot)b_{n_0}+\bar{z}(t)\bar{q}(\cdot)+W(z(t), \bar{z}(t), \theta)$ and denote by
$$F(0,U_t)|_{\mathscr{C}_0}=\tilde F(0,z,\bar{z}).$$
We write the Taylor formula
$$\tilde F(0,z,\bar{z})=\sum\limits_{i+j\geq2}\frac{\partial^{(i+j)}\tilde F}{\partial^i z\partial^j\bar z}\frac{z^i\bar z^j}{i!j!}.$$
Now the equation (\ref{zdot22}) on the center manifold is expressed by
\begin{equation}\label{zdot33}
\dot{z}(t)=\textrm{i}\omega_0 z(t)+g(z(t),\bar{z}(t)).\end{equation}
In order to calculate the normal form, we expand
\begin{equation}
g(z,\bar{z})=\sum\limits_{i+j\geq2}g_{i j}\frac{z^i\bar z^j}{i!j!}.\end{equation}

 According to the theory about normal form at a pair of imaginary roots \cite{Kuz}. We have the normal form till fifth order given by
 \begin{equation}
    \label{normalformfifth}
    \dot z=(\nu(\mu)+i)z+l_1(\mu)z^2\bar z+l_2(\mu)z^3{\bar z}^2+O(|z|^6)
 \end{equation}
with $\nu(\mu)=\frac{\alpha(\mu)}{\omega(\mu)}$,
\begin{equation}l_1(\mu)=\frac1{2\omega_0}\left({\rm Re} g_{21}-\frac{1}{\omega_0}{{\rm Im}}\{g_{20}g_{11}\}\right)\end{equation} and
\begin{equation}\begin{array}{rl}
&12l_2(\mu)=\frac 1 {\omega_0}{\rm Re}{g_{32}}\\&~~+\frac1{\omega_0^2}{{\rm Im}}\left\{g_{20}\bar g_{31}-g_{11}(4g_{31}+3\bar g_{22})-\frac13g_{02}(g_{40}+\bar g_{13})-g_{30}g_{12}\right\}\\
&~~+\frac{1}{\omega_0^3}{\rm Re}\left\{g_{20}\left[\bar g_{11}(3g_{12}-\bar g_{30})+g_{02}\left(\bar g_{12}-\frac13g_{30}\right)+\frac13\bar g_{02}g_{03}\right]\right\}\\
&~~+\frac{1}{\omega_0^3}{\rm Re}\left\{g_{11}\left[\bar g_{02}\left(\frac53\bar g_{30}+3g_{12}\right)+\frac13g_{02}\bar g_{03}-4g_{11}g_{30}\right]\right\}\\
&~~+\frac{3}{\omega_0^3}{{\rm Im}}\{g_{20}g_{11}\}{{\rm Im}}\{g_{21}\}+\frac1{\omega_0^4}{{\rm Im}}\{
g_{11}\bar g_{02}[\bar g_{20}^2-3\bar g_{20}g_{11}-4g_{11}^2\}\\
&~~+\frac1{\omega_0^4}{{\rm Im}}\left\{
g_{11} g_{20}\right\}[3{\rm Re}\{g_{11}g_{20}\}-2|g_{02}|^2].
 \end{array}\end{equation}
Notice that some  $g_{i j}$'s may depend on the form of $W_{i j}$. By using (\ref{abso}) and (\ref{uuuu}), we have
\begin{equation}\label{Wdot}\begin{array}{ll}\dot{W}&=\dot{U_t}-\dot{z}q b_{n_0}-\dot{\bar{z}}\bar{q}b_{n_0}\\&=\left\{\begin{array}{ll}
    A_0W-2{\rm Re}\{g(z,\bar{z})q(\theta)\}b_{n_0},& \theta\in[-r,0) \\
    A_0W-2{\rm Re}\{g(z,\bar{z})q(\theta)\}b_{n_0}+\tilde F,& \theta=0
\end{array}\right.\\&:= A_0W+H(z,\bar{z},\theta).
\end{array}
    \end{equation}
Again, expending $H(z,\bar z,\theta)$ as
\begin{equation}
H(z,\bar{z})=\sum\limits_{i+j\geq2}H_{i j}\frac{z^i\bar z^j}{i!j!},\end{equation}
we have
$$H_{i j}(\theta)=\left\{\begin{array}{ll}
-g_{i j}q(\theta)b_{n_0}-\bar{g}_{j i}\bar{q}(0)b_{n_0} ,& \theta\in[-r,0)\\
 -g_{i j}q(0)b_{n_0}-\bar{g}_{j i}\bar{q}(0)b_{n_0}+\frac{\partial^{(i+j)}\tilde F}{\partial^i z\partial^j \bar z},& \theta=0\end{array}
 \right.$$
and
 \begin{equation}\label{Waij}
 \begin{array}{l}
 (A_0-2\rm  i\omega_0 I)W_{20}(\theta)=-H_{20}(\theta),~A_0W_{11}(\theta)=-H_{11}(\theta),\\(A_0+2\rm  i\omega_0 I)W_{02}(\theta)=-H_{02}(\theta),\\
 (A_0-3\rm  i\omega_0 I)W_{30}(\theta)-3g_{20}W_{20}(\theta)-3\bar g_{02}W_{11}(\theta)=-H_{30}(\theta),\\
 (A_0-\rm  i\omega_0 I)W_{21}(\theta)-(g_{20}+2\bar g_{11})W_{11}(\theta)-2 g_{11}W_{20}(\theta)-\bar g_{02}W_{02}(\theta)=-H_{21}(\theta),\\
  (A_0+{\rm i}\omega_0 I)W_{12}(\theta)-(\bar g_{20}+2 g_{11})W_{11}(\theta)-2\bar g_{11}W_{02}(\theta)- g_{02}W_{20}(\theta)=-H_{12}(\theta),\\
   (A_0+3{\rm i}\omega_0 I)W_{03}(\theta)-3g_{02}W_{11}(\theta)-3\bar g_{20}W_{02}(\theta)=-H_{03}(\theta),\\
   (A_0-2{\rm i}\omega_0 I)W_{31}(\theta)-3(g_{20}+\bar g_{11})W_{21}(\theta)-(g_{30}+3\bar g_{12})W_{11}(\theta)-3  g_{11}W_{30}(\theta)\\~~~~-3g_{21}W_{20}(\theta)-\bar g_{03}W_{02}(\theta)-3\bar g_{02}W_{12}(\theta)=-H_{31}(\theta),\\
   A_0W_{22}(\theta)-g_{02}W_{30}(\theta)-(g_{20}+4\bar g_{11})W_{12}(\theta)-2(g_{21}+\bar g_{21})W_{11}(\theta)-2g_{12}W_{20}(\theta)
   \\~~~~-(\bar g_{20}+4g_{11})W_{21}(\theta)-2\bar g_{12}W_{02}(\theta)-\bar g_{02} W_{03}(\theta)=-H_{22}(\theta).
 \end{array}
 \end{equation}

Recalling the definition of $A_0$ in (\ref{inf2}), we know (\ref{Waij}) is actually a group of ordinary differential equations about $W_{ij}(\theta)$ with some given boundary conditions, solving which the explicit expressions of $W_{ij}$ and $g_{ij}$ are finally obtained.

\begin{figure}[htbp]
\begin{center}
a) \includegraphics[width=6cm]{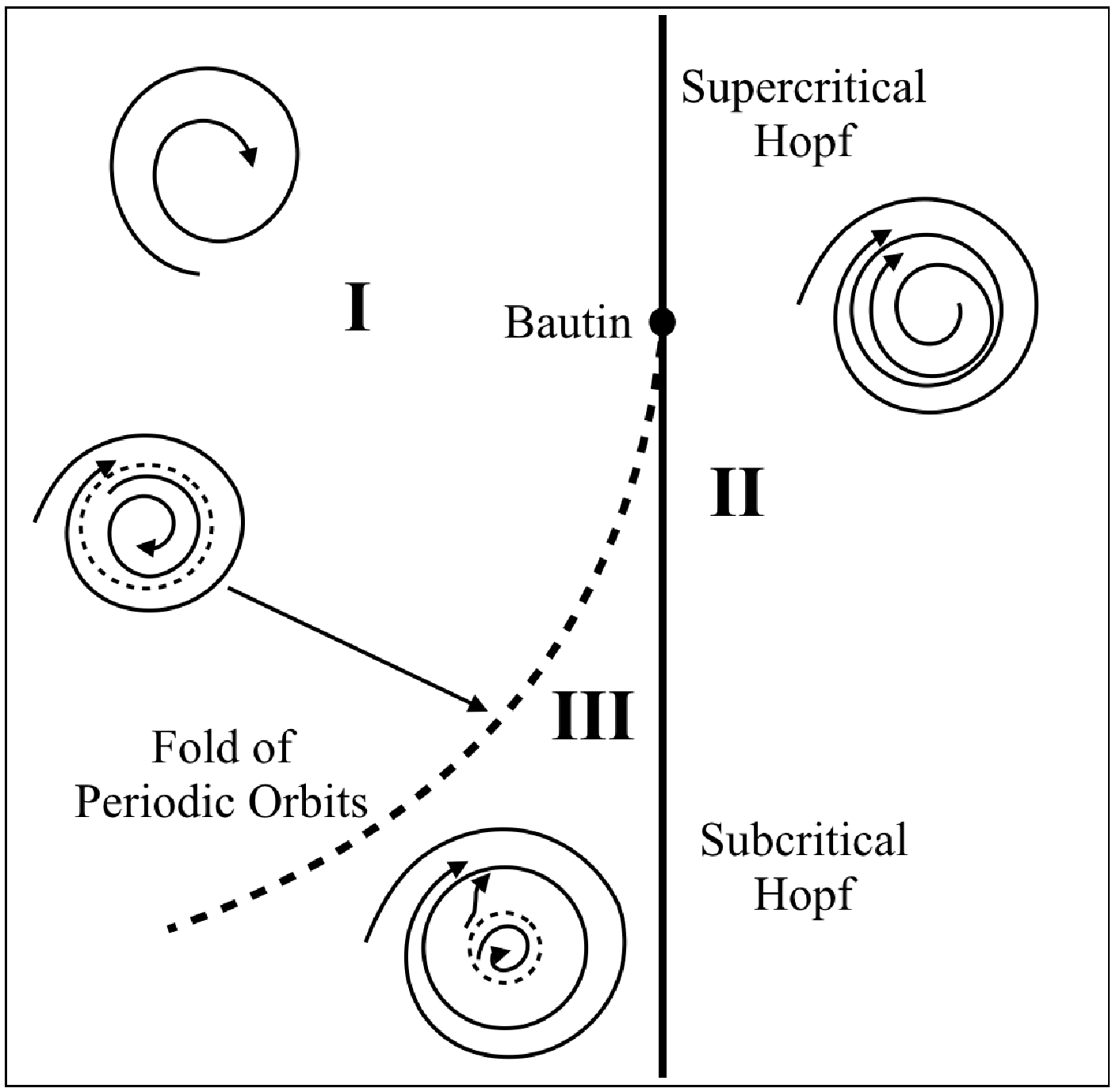} b) \includegraphics[width=6cm]{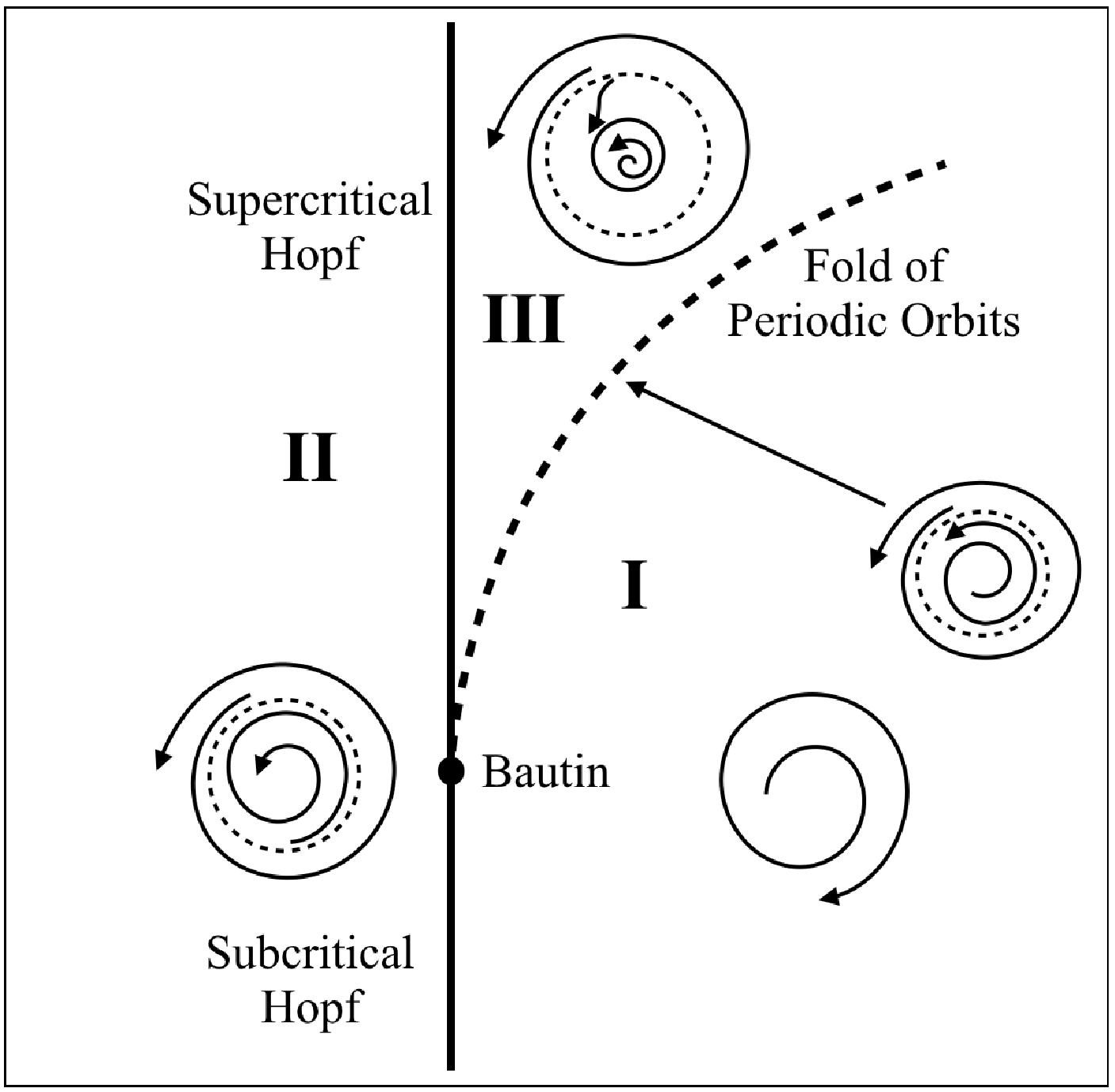} \end{center}
\caption{The Bautin bifurcation diagram near the equilibrium for a) $l_2(0)<0$ and b) $l_2(0)>0$. }\label{figbautin}
\end{figure}

According to the general results about Bautin bifurcations in \cite{Holmes}, we assume   $l_1(0)=0$, $l_2(0)\neq0$ and the map $\mu\mapsto(\nu(\mu),l_1(\mu))$ is regular at $\mu=0$, then the system undergoes a Bautin bifurcation. Precisely, the bifurcation diagram is given by Figure \ref{figbautin}.

Following the results given in \cite{Holmes}, we know that near the Bautin bifurcation, when $l_2(0)<0$ two limit cycles (region III) collide on the curve standing for fold bifurcation of periodic orbits, then disappear, and leave the system a stable equilibrium (region I), which undergoes a supercritical Hopf bifurcation and lead the system to stable periodic oscillations (region II). This is shown in Figure \ref{figbautin} (a). When $l_2(0)>0$, the situation is shown in Figure \ref{figbautin} (b), where the fold bifurcation of periodic orbits eliminates two limit cycles (region III), and leave the system an unstable equilibrium (region I), which undergoes a subcritical Hopf bifurcation and unstable periodic orbits appear (region II).

In the coming section, an example of Segel-Jackson model will be analyzed, where we find a Bautin bifurcation with $l_2(0)>0$ appears.

 \section {Example: Segel-Jackson model}
 In this section, we shall investigate the Bautin bifurcation for the Segel-Jackson model.
 \begin{equation}\label{S_J_model}
 \left\{
 \begin{array}{l}
 \frac{\partial u(x,t)}{\partial t}=d\Delta u(x,t)+[1+ku(x,t)]u(x,t)-au(x,t)v(x,t)\\
  \frac{\partial v(x,t)}{\partial t}=\Delta v(x,t)+u(x,t-\tau)v(x,t-\tau)-v^2(x,t)
 \end{array}~~~~~x\in[0,l\pi], t>0
 \right..
 \end{equation}
$u(x,t)$ stands for the density of the prey and $v(x,t)$ the predator. $d$ is the diffusion rate of the prey whereas the diffusion rate of predator is normalized by 1. $1+k u(x,t)$ stands for the reproduction rate per capita of prey. $a$ is the ability of the predator hunting the prey.

Equipping (\ref{S_J_model}) with homogeneous Neumann boundary condition $\left.\frac{\partial u}{\partial x}\right|_{x=0,l\pi}=\left.\frac{\partial v}{\partial x}\right|_{x=0,l\pi}=0$, there are two constant equilibria of system (\ref{S_J_model}): one is the trivial equilibrium $E_0=(0,0)$ and the other is the unique positive equilibrium $E_\ast=(u^\ast,v^\ast)=(\frac{1}{a-k},\frac{1}{a-k})$, provided that $k<a$.

The linear equation  of system (\ref{S_J_model}) at $E_\ast$ is
\begin{equation}\label{linear}
 \left\{
 \begin{array}{l}
 \frac{\partial u(x,t)}{\partial t}=d\Delta u(x,t)+\frac{k}{a-k}u(x,t)+\frac{a}{k-a}v(x,t),\\
  \frac{\partial v(x,t)}{\partial t}=\Delta v(x,t)+\frac{2}{k-a}v(x,t)+\frac{1}{a-k}u(x,t-\tau)+\frac{1}{a-k}v(x,t-\tau).
 \end{array}
 \right.
 \end{equation}
Then we have the
  characteristic equation of system (\ref{linear}) is
\begin{equation}\label{characteristic}
 \lambda^2+A_n\lambda+B\lambda{{\rm e}}^{-\lambda\tau}+C_n+D_n{{\rm e}}^{-\lambda\tau}=0,
 \end{equation}
 with $A_n=\frac{k-2}{k-a}+\frac{(d+1)n^2}{l^2}$, $B=\frac{1}{k-a}$, $C_n=\left(\frac{k}{a-k}-\frac{d n^2}{l^2}\right)\left(\frac{2}{k-a}-\frac{n^2}{l^2}\right)$, $D_n=\frac{k+a}{(a-k)^2}-\frac{d n^2}{(a-k) l^2 }$.

 In the following, we will discuss the stability of $E_\ast$, regarding time delay $\tau$ and $k$ as  bifurcation parameters.

 When $\tau=0$, Eq.(\ref{characteristic}) becomes
 \begin{equation}\label{tau0}
 \lambda^2+(A_n+B)\lambda+C_n+D_n=0.
 \end{equation}

 In fact, this is just the case investigated in \cite{Wang}, and we state the main results here: when $\tau=0$, assume that  $$k<\min\{1,a,d\}~~~~~~(H_1)$$
 holds true, then the positive equilibrium $E_\ast$ is locally asymptotically  stable.

In the following part, we will investigate the existence of Hopf bifurcations destabilizing $E_\ast$ as $\tau$ increases, where we always assume $(H_1)$ holds true.

When $\tau>0$, following the method given in \cite{wei2}, we plug $\lambda={\rm i}\omega$ into Eq.(\ref{characteristic}) and separate the real and imaginary parts, then obtain
\begin{equation}\label{fl}
\left\{\begin{array}{l}
D_n\cos\omega\tau+B \omega\sin\omega\tau=\omega^2-C_n,\\
-D_n\sin\omega\tau+B \omega\cos\omega\tau=-A_n\omega.
\end{array}
\right.
\end{equation}
This is a linear equation about unknowns $\cos\omega\tau$ and $\sin\omega\tau $. Once $\omega$ and $n$ are fixed, we can solve from (\ref{fl}) that $\cos\omega\tau:=\tilde C_n(\omega)$ and $\sin\omega\tau:=\tilde S_n(\omega).$

 Adding the square of both sides of Eq.(\ref{fl}), we have an algebraic equation about $\omega^2$, denoted by
 \begin{equation}\label{omega}
(\omega^2)^2+M_1\omega^2+M_2=0
\end{equation}
with $M_1=-2 C_n+A_n^2-B^2$ and $M_2=C_n^2-D_n^2$.

If
$M_1<0,~ M_2>0, ~M_1^2-4 M_2>0$
hold true,
Eq.(\ref{omega}) has two positive roots, given by $$\omega=\omega_n^\pm=\sqrt{\frac{-M_1\pm\sqrt{M_1^2-4M_2}}{2}}.$$
Else if $M_2<0$ holds true, Eq.(\ref{omega}) has only one positive root,  denoted by
$$\omega=\omega_n^+=\sqrt{\frac{-M_1+\sqrt{M_1^2-4M_2}}{2}}.$$

Plugging $\omega_n^+$ or $\omega_n^\pm$ into Eq.(\ref{fl}), at most two sequences of possible  Hopf bifurcation values are obtained, which are denoted by
\begin{equation}\label{tauasta}
    \tau_{n,j}^\pm=\left\{ \begin{array}{ll}
    \frac{\arccos \tilde C_n(\omega_n^\pm)+2j\pi}{\omega_n^\pm},&~if~ \tilde  S_n(\omega_n^\pm)\geq0\\
    \frac{\arcsin \tilde S_n(\omega_n^\pm)+2\pi+2j\pi}{\omega_n^\pm},& ~if~ \tilde S_n(\omega_n^\pm)<0 ~and~ \tilde C_n(\omega_n^\pm) \geq0\\
    \frac{-\arcsin \tilde S_n(\omega_n^\pm)+\pi+2j\pi}{\omega_n^\pm}, &~if~ \tilde S_n(\omega_n^\pm)<0 ~and~ \tilde C_n(\omega_n^\pm)<0
    \end{array} \right. ~~~~n,j=0,1,2,\cdots.
\end{equation}

We know when $\tau=\tau_{n,j}^\pm, n,j=0,1,2,\cdots$, Eq.(\ref{characteristic}) has a pair of imaginary roots.
By re-scaling time $t\rightarrow t/\tau$, making transformation
$\tilde u=u(x,t)-u^\ast,$ $\tilde v=v(x,t)-v^\ast$
and dropping the tildes for convenience, we have system (\ref{S_J_model}) becomes
\begin{equation}\label{rescaling}
\left\{\begin{array}{l}
\frac{\partial u}{\partial t}=\tau\left[d\Delta u+(1+2 k u^\ast-a v^\ast) u-a u^\ast v+f_1(u_t,v_t)\right]\\
\frac{\partial v}{\partial t}=\tau\left[\Delta v-2 v^\ast v+v^\ast u_t(-1) +u^\ast v_t(-1)+f_2(u_t,v_t)\right]
\end{array}
\right.
\end{equation}
where
$u=u(x,t)$, $v=v(x,t)$.  For any $(\phi_1,\phi_2)\in\mathscr{C}:=C([-1,0],X)$ and $X:=\{(u,v)\in H^2(0,l\pi)\times H^2(0,l\pi):\frac{du}{dx}=\frac{dv}{dx}=0~ at~ x=0,l\pi\}$, two nonlinear functionals are defined by
\begin{equation}\label{f1}f_1(\phi_1,\phi_2)=k{\phi_1}^2(0)-a{\phi_1}(0){\phi_2}(0)\end{equation}
and
\begin{equation}\label{f2}f_2(\phi_1,\phi_2)=-\phi_2^2(0)+\phi_1(-1)\phi_2(-1).\end{equation}

We will calculate the first and second Lyapunov coefficients of the normal forms when $\tau=\tau^\ast\in\{\tau_{n,j}^\pm,~n,j=0,1,2,\cdots\}$ and for any fixed $k$, while Eq.(\ref{characteristic})
has a pair of imaginary roots $\lambda=\pm\textrm{i}\omega_n^\pm$ for $n=n_0$. For simplification, we denote  $\omega_{n_0}^\pm$ by $\omega^\ast$.

Let $\tau=\tau^\ast+\epsilon$. We can rewrite system (\ref{rescaling}) in an abstract form in the phase space $\mathscr{C}:=C([-1,0],X)$ as

 \begin{equation}\label{abs}\dot U(t)=\tilde D\Delta U(t)+L_\epsilon(U_t)+F(\epsilon,U_t)\end{equation}
  where~
 $\tilde D=(\tau^\ast+\epsilon)D$, $D=\left(\begin{array}{cc}
      d & 0 \\
      0 & 1
 \end{array}\right)$ and $L_\epsilon:\mathscr{C}\rightarrow X$, $F:\mathscr{C}\rightarrow X$ are defined, respectively, by
 $$L_\epsilon(\phi)=(\tau^\ast+\epsilon)B_1\phi(0)+(\tau^\ast+\epsilon)B_2\phi(-1)$$
 and
 $$F(\epsilon,\phi)=(F_1(\epsilon,\phi(\theta)),F_2(\epsilon,\phi(\theta)))^{\textrm{T}}$$
with
 $$B_1=\left(\begin{array}{cc}
 1+2 k u^\ast-a v^\ast& -a u^\ast \\
      0 & -2 v^\ast\end{array}\right),~~B_2=\left(\begin{array}{cc}
 0& 0 \\
      v^\ast& u^\ast\end{array}\right),$$
\begin{equation}\label{F1}
   F_1(\epsilon,\phi)=(\tau^\ast+\epsilon)f_1(\phi_1(\theta),\phi_2(\theta)),~~~for~any~ \phi=(\phi_1,\phi_2)\in\mathscr{C}
\end{equation}
and
\begin{equation}\label{F2}
    F_2(\epsilon,\phi)=(\tau^\ast+\epsilon)f_2(\phi_1(\theta),\phi_2(\theta)),~~~for~any~ \phi=(\phi_1,\phi_2)\in\mathscr{C}.
\end{equation}
Thus, the linearized equation of system (\ref{rescaling}) at the equilibrium $(0,0)$ has the form
 \begin{equation}\label{rescalingl} \dot U(t)=\tilde D\Delta U(t)+L_\epsilon(U_t).  \end{equation}

 We have that the solution operator of system (\ref{rescaling}) forms a $C_0$ semigroup with the infinitesimal generator $A_\epsilon$,

 \begin{equation}\label{inf} A_\epsilon \phi=\left\{\begin{array}{ll} \dot\phi(\theta),&\theta\in [-r,0),\\
 \tilde D\Delta\phi(0)+L_\epsilon(\phi),&\theta=0.\end{array}\right.\end{equation}
 The domain of $A_\epsilon$ is chosen by

 $${\rm dom}(A_\epsilon):=\{\phi\in\mathscr{C}:\dot{\phi}\in\mathscr{C}, \phi(0)\in {\rm dom}(\Delta), \dot{\phi}(0)= \tilde D\Delta\phi(0)+L_\epsilon(\phi)\}.$$

In order to use the results we obtained in the precious section, in the space $\mathscr{C}$ we rewrite Eq.(\ref{abs}) as the abstract form
 \begin{equation}\label{abso2}\dot{U}(t)=A_\epsilon U_t+R(\epsilon,U_t)\end{equation}
where the nonlinear term is  given by
 $$R(\epsilon,U_t)(\theta)=\left\{\begin{array}{ll}
      0,&   \theta\in[-1,0), \\
      F(\epsilon,U_t),& \theta=0.
 \end{array}\right.$$

 In the state space $X$, we know, $\left(\cos\frac{n x}{l}(0,1)^\textrm{T},\cos\frac{n x}{l}(1,0)^\textrm{T}\right),~~n=0,1,2,\cdots,$ are eigenfunctions of $\Delta$
 with the no-flux boundary conditions. Moreover, they form a basis of $X$ which are normalized by
 $b_n=\frac{\cos(\frac{nx}{l})}{\|\cos\frac{nx}{l}\|_{L_2}}$.
 For convenience we use $\beta_n=\{ b_n(1,0)^{\textrm{T}},b_n(0,1)^{\textrm{T}}\}$ to be the basis of $X$.

Using the notations introduced in Section 2, for any $\phi=(\phi^{(1)},\phi^{(2)})^{\textrm{T}}\in\mathscr{C}$, we denote
 $$\phi_n=<\phi,\beta_n>=\left(<\phi^{(1)},b_n>,<\phi^{(2)},b_n>\right)^{\textrm{T}}$$
 as the coordinates of $\phi$ on $\beta_n$ in $X$.

 The restriction of $A$, $A_{\epsilon,n}$ is then defined by
 $$A_{\epsilon,n}(\phi_n(\theta)b_n)=\left\{ \begin{array}{ll}
      \dot\phi_n(\theta)b_n,& \theta\in[-1,0), \\
      \int_{-1}^0 d\eta_n(\epsilon,\theta)\phi_n(\theta)b_n,& \theta=0,
 \end{array}\right.$$
with
 $$\eta_n(\epsilon,\theta)=\left\{ \begin{array}{ll}
 -(\tau^\ast+\epsilon)B_2,& \theta=-1,\\
 0,& \theta\in(-1,0),\\
 (\tau^\ast+\epsilon)(B_1-\frac{n^2}{l^2}D),& \theta=0.\end{array}\right.$$

Similarly, the linear operator $L_2$ has the restriction
 $$L_{\epsilon,n}(\phi_n)=(\tau^\ast+\epsilon)B_1\phi_n(0)+(\tau^\ast+\epsilon)B_2\phi_n(-1).$$
 Obviously, we have, from the linear Eq.(\ref{rescalingl})
 $$-\frac{n^2}{l^2}\tilde D\phi_n(0)+L_{\epsilon,n}(\phi_n)=\int_{-1}^0 d\eta_n(\epsilon,\theta)\phi_n(\theta)b_n.$$

Similar as the previous section we can define the adjoint operator $A^\ast$ of $A_0$ as

 $$A^\ast\varphi(s)=\left\{\begin{array}{ll}
 -\dot\varphi(s),& s\in(0,1],\\
 \sum\limits_{n=0}^\infty\int_{-1}^0d\eta_n^\textrm{T}(0,t)\varphi_n(-t)b_n,& s=0.\end{array}
 \right.$$
As we stated at the beginning of this section, we can calculate the eigenfunctions of $A_0$ and $A^\ast$ at the eigenvalue $\textrm{i}\omega^\ast\tau^\ast$ and $-\textrm{i}\omega^\ast\tau^\ast$, which are, respectively, $qb_{n_0}$ and $q^\ast b_{n_0}$.

By direct calculations, we have
$q(\theta)=(1,q_1)^\textrm{T}\textrm{e}^{\textrm{i}\omega^\ast\tau^\ast\theta}$ and   $q^\ast(s)=M(q_2,1)^\textrm{T}\textrm{e}^{\textrm{i}\omega^\ast\tau^\ast s}$.
Here
$$q_1=\frac{1+2 k u^\ast-a v^\ast-\frac{d n^2}{l^2}+\textrm{i}\omega^\ast\tau^\ast}{a u^\ast},$$  $$q_2=\frac{-2 v^\ast-\frac{n^2}{l^2}+\textrm{e}^{\textrm{i}\omega^\ast\tau^\ast}u^\ast+\textrm{i}\omega^\ast\tau^\ast}{a u^\ast}$$
and
$$\bar{M}=[{(\bar{q_2}+q_1)+\textrm{e}^{-\textrm{i}\omega^\ast\tau^\ast}\tau^\ast(v^\ast+q_1 u^\ast)}]^{-1},$$
such that $(q^\ast,q)_c=1$.  By the formal adjoint theory, we have
$(\bar{q^\ast},q)_c=(q^\ast,\bar{q})=0.$ By using the standard theory of phase space decomposition, we have the eigenvalues $\Lambda=\{\pm\textrm{i}\omega^\ast\tau^\ast\}$ can be used to decompose $\mathscr{C}$ by $\mathscr{C}=P\bigoplus Q$. The center space $P$ is given by
$$P=\{z q b_{n_0}+\bar z q b_{n_0}|z\in \mathbb C\}$$
and its orthogonal complement space is
$$Q=\{\phi\in \mathscr{C}|(q^\ast b_{n_0},\phi)=0 ~and ~(\bar{q}^\ast b_{n_0},\phi)=0\}.$$
Thus, the state of system (\ref{abso2}) could be decomposed, correspondingly, by
$$U_t=z(t)q(\cdot)b_{n_0}+\bar{z}(t)\bar{q}(\cdot)b_{n_0}+W(t,\cdot)$$
where $W(t,\cdot)\in Q.$

Notice the relation between $q^\ast$ and $q$, we calculate $z(t)$ by
\begin{equation}\label{z}z(t)=(q^\ast b_{n_0},U_t).\end{equation}
Thus, $W(t,\theta)=U_t(\theta)-2{\rm Re}\{z(t)q(\theta)b_{n_0}\}$.
Then we obtain the reduced equation onto the center space
\begin{equation}\label{zdot2}\dot{z}(t)=\textrm{i}\omega_0z(t)+\bar{q}^{\ast\textrm{T} }<F(0,U_t),\beta_{n_0}>,
\end{equation}
where the nonlinear term is given by
$$<F,\beta_n>:=(<F_1,b_n>,<F_2,b_n>)^\textrm{T},$$
with $F_1$ and $F_2$ already defined in (\ref{F1}) and (\ref{F2}).

Using the center manifold theorem \cite{Lin}, there exists a center manifold $\mathscr{C}_0$ and on $\mathscr{C}_0$ we have $W(t,\cdot)$ can be expressed by series
\begin{equation}\label{cm}\begin{array}{ll}W(t,\cdot)&=W(z(t),\bar{z}(t),\cdot)=W_{20}\frac{z^2}{2}+W_{11} z\bar{z}+W_{02}\frac{\bar{z}^2}{2}+W_{30}\frac{z^3}{6}+W_{21}\frac{z^2\bar{z}}{2}+W_{12}\frac{z\bar{z}^2}{2}\\&+W_{03}\frac{\bar{z}^3}{6}+W_{31}\frac{z^3\bar{z}}{6}+W_{22}\frac{z^2\bar{z}^2}{4}+\cdots.\end{array}
\end{equation}
 For solution $U_t\in\mathscr{C}_0,$ we have $U_t=z(t)q(\cdot)b_{n_0}+\bar{z}(t)\bar{q}(\cdot)+W(z(t), \bar{z}(t), \theta)$ and denote by
$$F(0,U_t)|_{\mathscr{C}_0}=\tilde F(0,z,\bar{z}).$$
We write the Taylor expression
\begin{equation}\label{F}\begin{array}{ll}
    \tilde F(0,z,\bar{z})&=\tilde F''_{z z}\frac{z^2}{2}+\tilde F''_{z\bar{z}}z\bar{z}+\tilde F''_{\bar{z}\bar{z}}\frac{\bar{z}^2}{2}+\tilde F'''_{z z z}\frac{z^3}{6}+\tilde F'''_{z z \bar z}\frac{z^2\bar z}{2}+\tilde F'''_{z \bar z \bar{z}}\frac{z\bar{z}^2}{2}+\tilde F'''_{\bar{z} \bar{z}\bar{z}}\frac{\bar{z}^3}{6}\\&+\tilde F''''_{z z z z}\frac{z^4}{24}+\tilde F''''_{z z z\bar{z}}\frac{z^3\bar{z}}{6}+\tilde F''''_{z z\bar{z}\bar{z}}\frac{z^2\bar{z}^2}{4}+\tilde F''''_{z \bar z\bar{z}\bar z}\frac{z\bar{z}^3}{6}+\tilde F''''_{\bar{z}\bar{z}\bar{z}\bar{z}}\frac{\bar{z}^4}{24}\\&+\tilde F'''''_{z z z\bar{z}\bar{z}}\frac{z^3\bar{z}^2}{12}+\cdots.\end{array}\end{equation}
In fact, from  Eqs.(\ref{f1}-\ref{f2})   and (\ref{F1}-\ref{F2}), we have
\begin{equation}\label{FFF}
\begin{array}{l}
\tilde{F}''_{zz}=\tau^\ast\left(\begin{array}{c}2k-2aq_1\\
-2q_1^2+2q_1{{\rm e}}^{-2{{\rm i}}\omega^\ast\tau^\ast}\end{array}\right)b_{n_0}^2,\\
\tilde{F}''_{\bar z\bar z}=\tau^\ast\left(\begin{array}{c}2k-2a\bar q_1\\
-2\bar q_1^2+2\bar q_1{{\rm e}}^{2{{\rm i}}\omega^\ast\tau^\ast}\end{array}\right)b_{n_0}^2,\\
\tilde{F}''_{  z\bar z}=\tau^\ast\left(\begin{array}{c}2k-a(q_1+\bar q_1)\\
-2\bar q_1 q_1+(q_1+\bar q_1)\end{array}\right)b_{n_0}^2,
\end{array}
\end{equation}
and all the rest higher-order derivatives are zero.

Now rewrite the equation on the center manifold (\ref{zdot2}) as
\begin{equation}\label{zdot3}
\dot{z}(t)=\textrm{i}\omega_0z(t)+g(z,\bar{z}).\end{equation}
In order to calculate the normal form, we expand
\begin{equation}\label{g}\begin{array}{ll}
g(z,\bar{z})&=g_{20}\frac{z^2}{2}+g_{11}z\bar{z}+g_{02}\frac{\bar{z}^2}{2}+g_{21}\frac{z^2\bar{z}}{2}+g_{12}\frac{z\bar{z}^2}{2}+g_{03}\frac{\bar{z}^3}{6}+g_{40}\frac{z^4}{24}+g_{31}\frac{z^3\bar{z}}{6}+g_{22}\frac{z^2\bar{z}^2}{4}\\&+g_{13}\frac{z\bar{z}^3}{6}+g_{04}\frac{\bar{z}^4}{24}+g_{32}\frac{z^3\bar{z}^2}{12}+\cdots.\end{array}\end{equation}
By comparing the same order of terms, we obtain the expressions of $g_{ij}$'s, which, together with the detailed derivations, are given in  the Appendix.

According to the results given in Section 2, the first and second Lyapunov coefficients can be expressed by
\begin{equation}\label{l1l1}l_1(k,\tau^\ast)=\frac1{2\omega^\ast\tau^\ast}\left({\rm Re} g_{21}-\frac{1}{\omega^\ast\tau^\ast}{{\rm Im}}\{g_{20}g_{11}\}\right)\end{equation}
and
\begin{equation}\label{l2l2}\begin{array}{rl}
&12l_2(k,\tau^\ast)=\frac 1 {\tau^\ast\omega^\ast}{\rm Re}{g_{32}}\\&~~+\frac1{(\tau^\ast\omega^\ast)^2}{{\rm Im}}\left\{g_{20}\bar g_{31}-g_{11}(4g_{31}+3\bar g_{22})-\frac13g_{02}(g_{40}+\bar g_{13})-g_{30}g_{12}\right\}\\
&~~+\frac{1}{(\tau^\ast\omega^\ast)^3}{\rm Re}\left\{g_{20}\left[\bar g_{11}(3g_{12}-\bar g_{30})+g_{02}\left(\bar g_{12}-\frac13g_{30}\right)+\frac13\bar g_{02}g_{03}\right]\right\}\\
&~~+\frac{1}{(\tau^\ast\omega^\ast)^3}{\rm Re}\left\{g_{11}\left[\bar g_{02}\left(\frac53\bar g_{30}+3g_{12}\right)+\frac13g_{02}\bar g_{03}-4g_{11}g_{30}\right]\right\}\\
&~~+\frac{3}{(\tau^\ast\omega^\ast)^3}{{\rm Im}}\{g_{20}g_{11}\}{{\rm Im}}\{g_{21}\}+\frac1{(\tau^\ast\omega^\ast)^4}{{\rm Im}}\{
g_{11}\bar g_{02}[\bar g_{20}^2-3\bar g_{20}g_{11}-4g_{11}^2\}\\
&~~+\frac1{(\tau^\ast\omega^\ast)^4}{{\rm Im}}\left\{
g_{11} g_{20}\right\}[3{\rm Re}\{g_{11}g_{20}\}-2|g_{02}|^2].
 \end{array}\end{equation}

\begin{figure}[htbp]
\begin{center}
 a)\includegraphics[width=7cm]{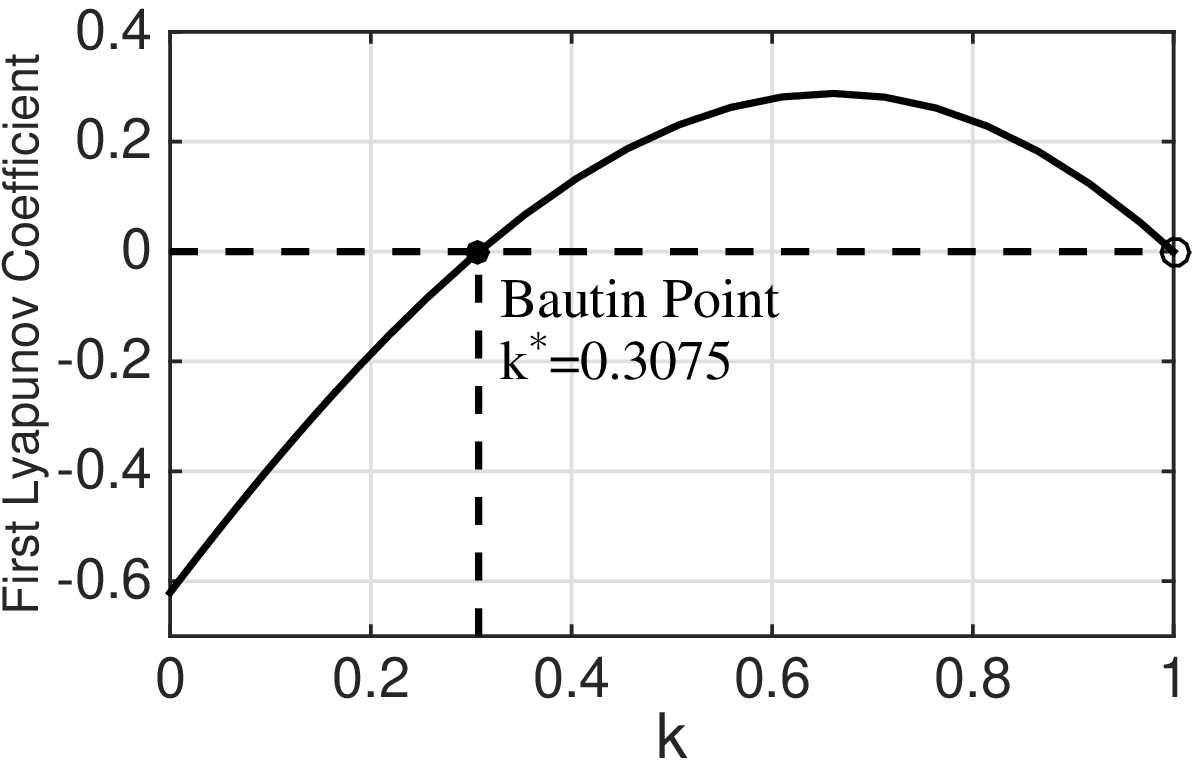} b)\includegraphics[width=7cm]{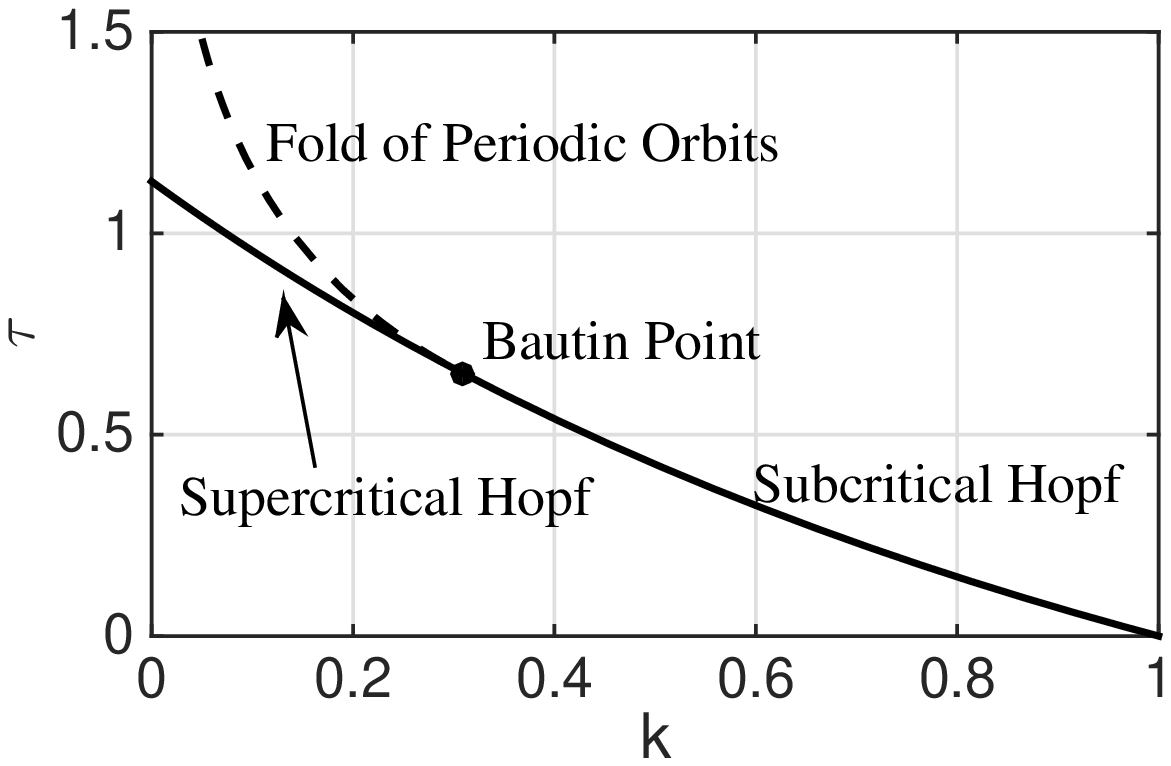} \end{center}
\caption{ $a= 5$. a) the first Lyapunov coefficient $l_1$ varies with respect to $k$. b) the Bautin bifurcation diagram near the Bautin point $(k^\ast, \tau^\ast)=(0.3075, 0.6543)$.}\label{fig1}
\end{figure}

Fixing $a=5$ and $d=1$, we can use the method above to determine the Hopf bifurcation values of $\tau$ for any given $k$. By using the formula (\ref{tauasta}) , we know at $(k,\tau)=(k^\ast,\tau^\ast)=(0.3075,0.6543)$,   the characteristic equation (\ref{characteristic}) has a pair of imaginary roots $\lambda=\pm 0.4233{{\rm i}}$, i.e., $\omega^\ast=0.4233$, and all the rest roots have negative real part. Further calculations from Eq.(\ref{l1l1})  and (\ref{l2l2}) yield $l_1(k^\ast,\tau^\ast)=0$ (as shown in Figure \ref{fig1} (a)) and  $l_2(k^\ast,\tau^\ast)=0.0376>0$.

\begin{figure}[htbp]
\begin{center}
 \includegraphics[width=7cm]{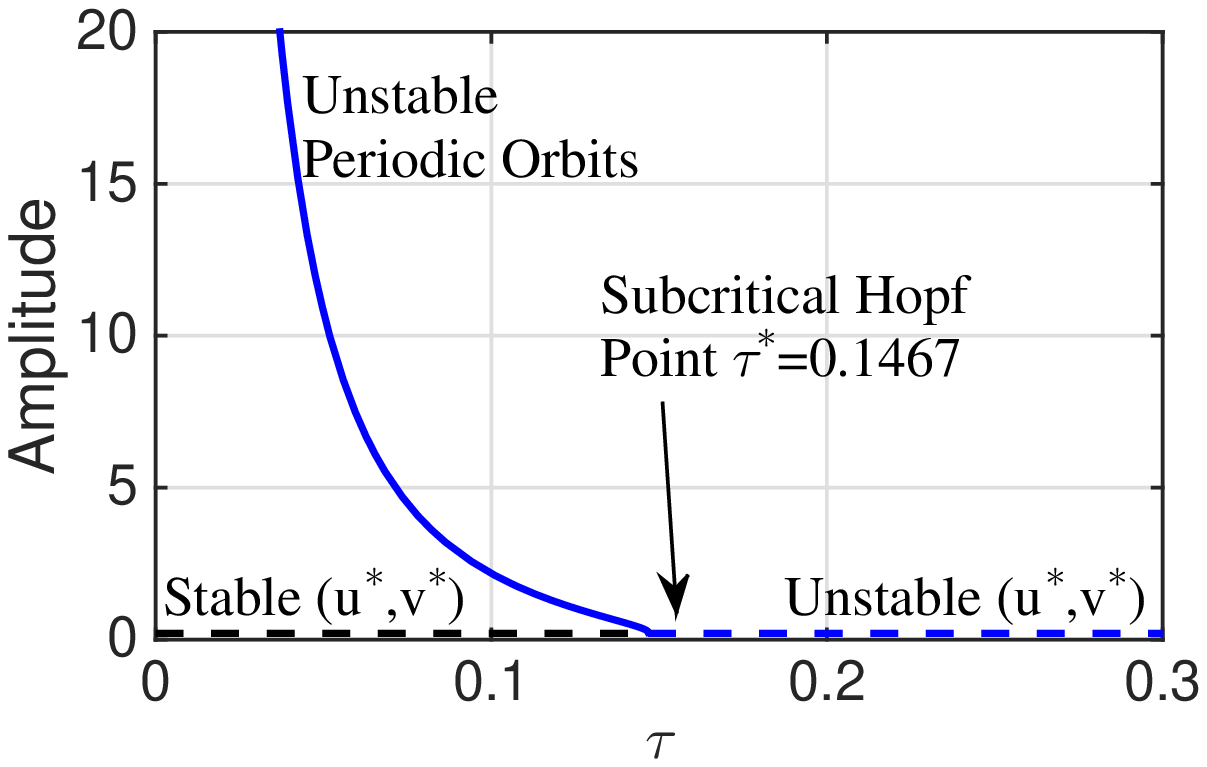} \includegraphics[width=7cm]{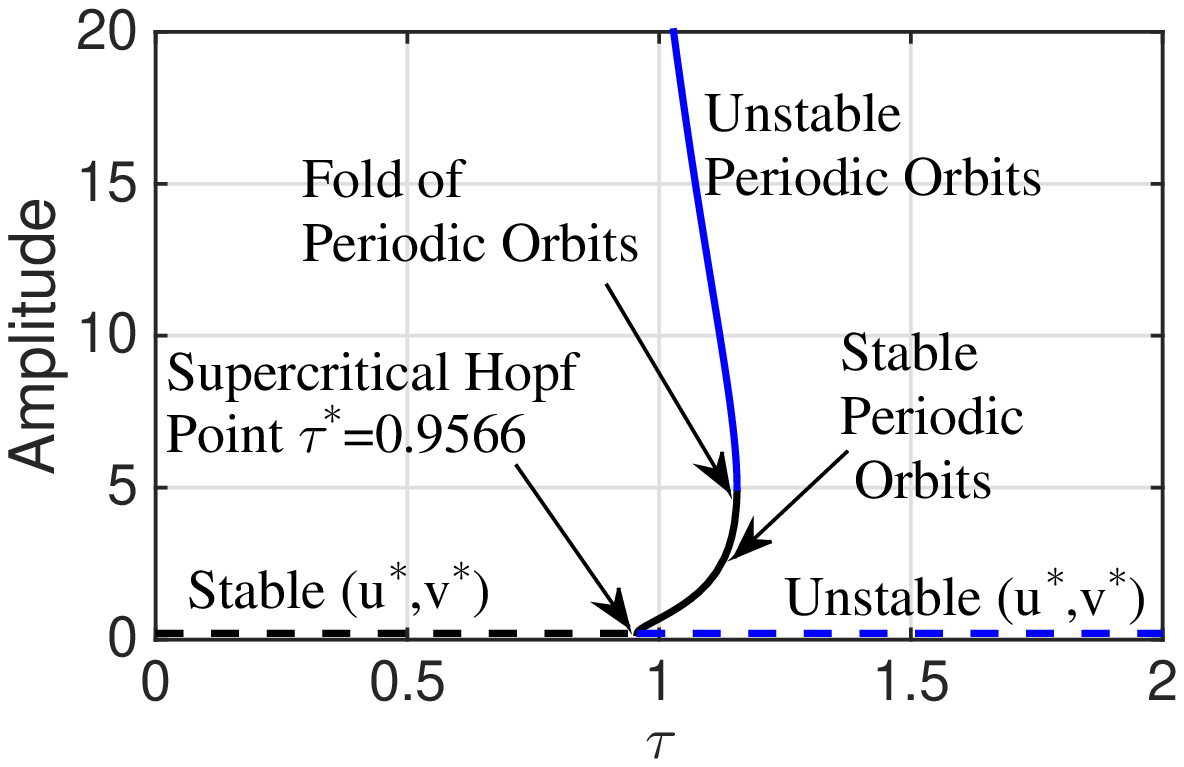} \end{center}
\caption{$a= 5$. (a)  Subcritical bifurcation diagram when $k=0.8$. (b) Supercritical bifurcation diagram when $k=0.1$.}\label{fig4}
\end{figure}

The transversality condition can be verified, numerically, by
\begin{equation}\label{verytran}
\left|\begin{array}{cc}
\frac{\partial\frac{\alpha(k,\tau)}{\omega(k,\tau)}}{\partial k} & \frac{\partial\frac{\alpha(k,\tau)}{\omega(k,\tau)}}{\partial \tau}\\
\frac{\partial l_1(k,\tau)}{\partial k} &\frac{\partial l_1(k,\tau)}{\partial  \tau}
\end{array}\right|
_{(k,\tau)={k^\ast,\tau^\ast}}=\left|\begin{array}{cc}
 \frac{\alpha'_k(k^\ast,\tau^\ast) }{\omega(k^\ast,\tau^\ast)}& \frac{\alpha'_\tau(k^\ast,\tau^\ast) }{\omega(k^\ast,\tau^\ast)} \\
\frac{\partial l_1(k,\tau)}{\partial k} &\frac{\partial l_1(k,\tau)}{\partial \tau}
\end{array}\right|=\left|\begin{array}{cc}0.2946&0.3085\\
1.5728&-1.1954\end{array}\right|=-0.8374\neq0.\end{equation}
 Thus the map $(k,\tau)\mapsto(\frac{\alpha(k,\tau)}{\omega(k\tau)},l_1(k,\tau))$ is regular at $(k^\ast,\tau^\ast)$, and  the system undergoes a Bautin bifurcation at $(k^\ast,\tau^\ast)$ with  the bifurcation diagram topologically equivalent to Figure \ref{figbautin} (b). When $k<k^\ast=0.3075$, the Hopf bifurcation is subcritical, as shown in Figure \ref{fig4} (a). The Bautin bifurcation diagrams is shown in Figure \ref{fig1} (b). When $k>k^\ast$, the Hopf bifurcation is supercritical. As shown in Figure \ref{fig4} (b), the fold bifurcation of periodic orbits eliminates the stable and unstable periodic orbits.

\begin{figure}[htbp]
\begin{center}
 \includegraphics[width=7cm]{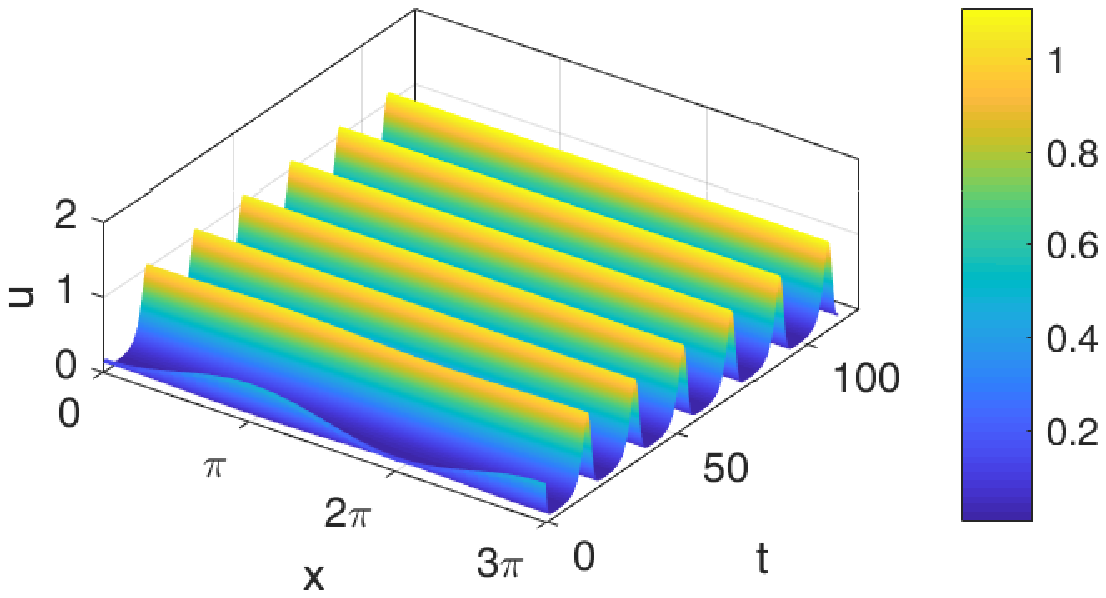} \includegraphics[width=7cm]{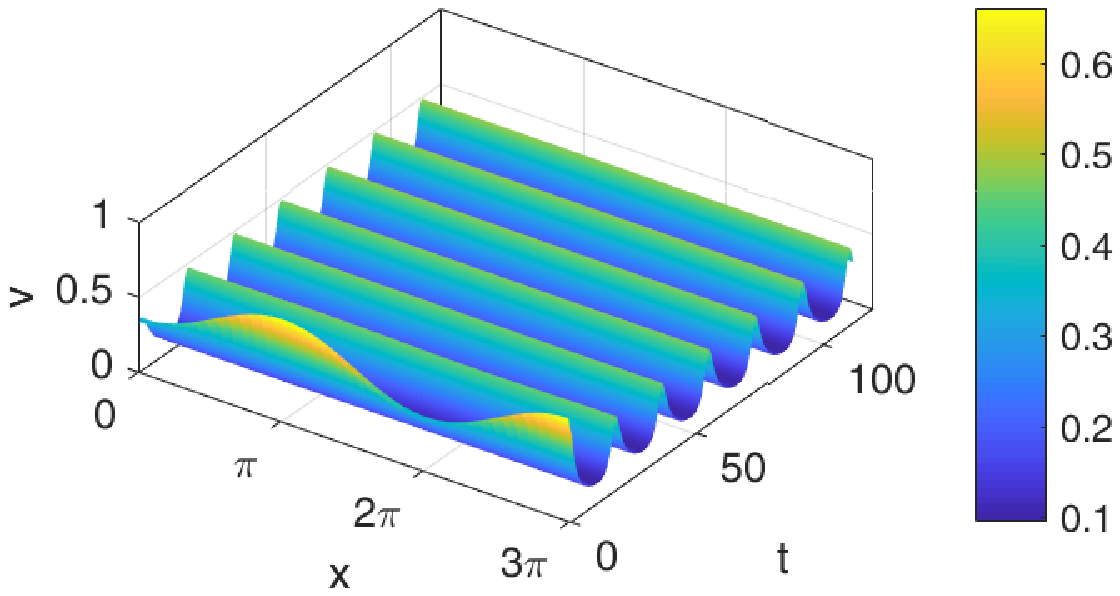} \end{center}
\caption{ $k=0.1$ $a= 5$. For $\tau=1.05$, the solution of (\ref{S_J_model}) with initial value  $u_0(x,t)=0.3-0.16\cos x$,  $v_0(x,t)=0.5-0.16\cos x$  converges to a periodic solution.}\label{fig6}
\end{figure}

\begin{figure}[htbp]
\begin{center}
 \includegraphics[width=7cm]{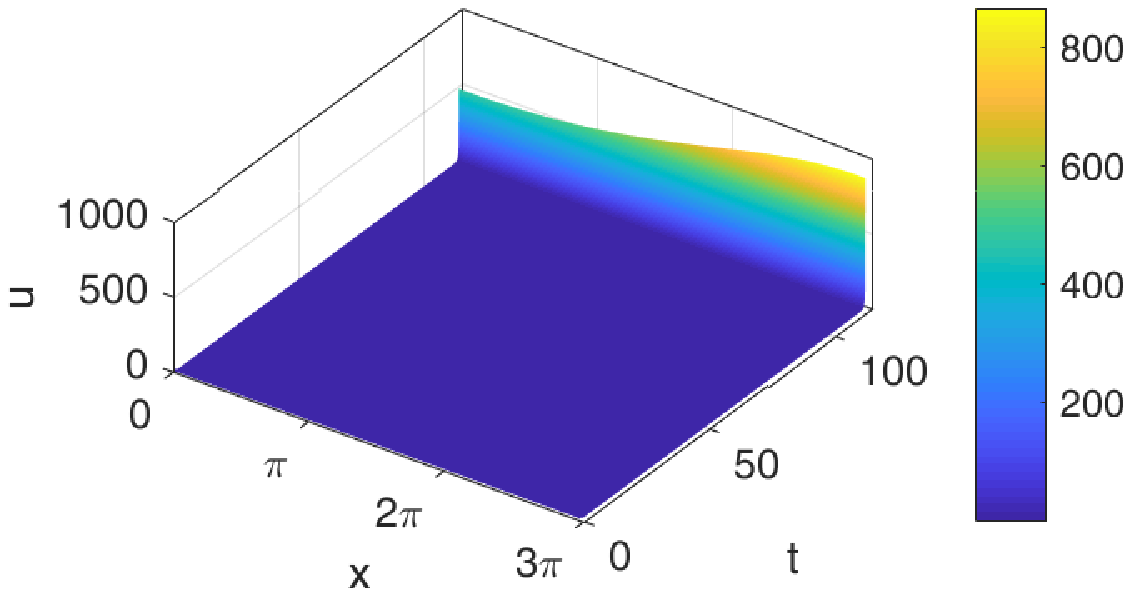} \includegraphics[width=7cm]{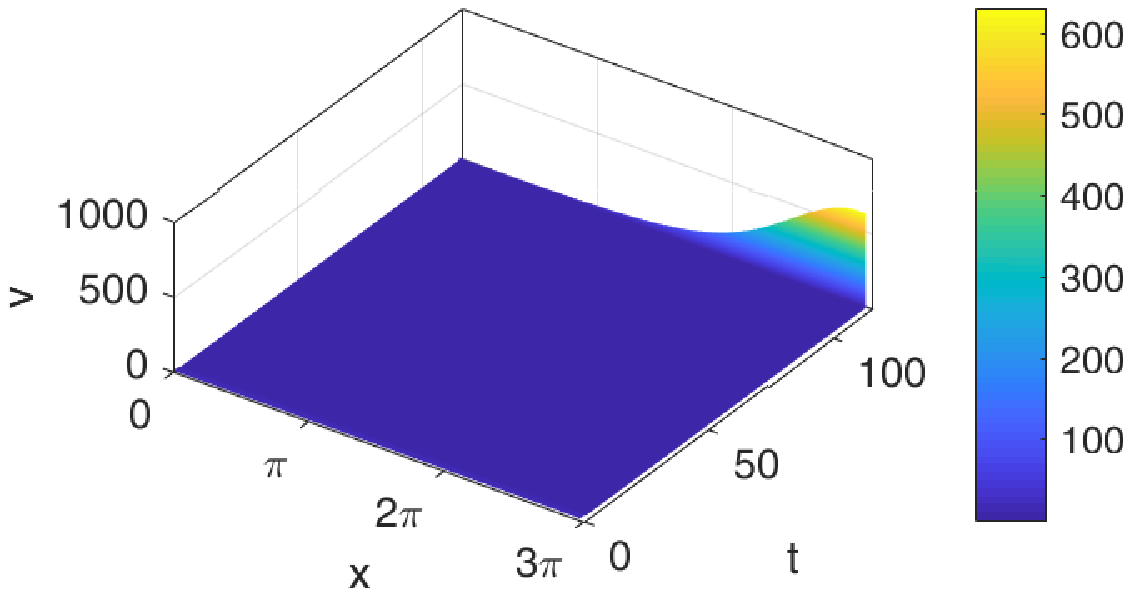} \end{center}
\caption{ $k=0.1$, $a= 5$. For $\tau=1.05$, the solution of (\ref{S_J_model}) with initial value   $u_0(x,t)=10.3-0.16\cos x$,  $v_0(x,t)=10.5-0.16\cos x$  diverges to infinity.}\label{fig7}
\end{figure}

 In fact, we find that when $k=0.8$ and $\tau\in(0.9566,1.1550)$, solutions with small initial values converge to a stable periodic orbits (Figure \ref{fig6}), and those with large initial values diverge to infinity as time goes to infinity (Figure \ref{fig7}).

\section*{Conclusion} In this paper, a universal and explicit method to calculate the first and second Lyapunov coefficients at a pair of imaginary eigenvalues in reaction-diffusion system with time delays is given, which
can be used to determine the dynamics near a Bautin bifurcation point. As an example, the method is applied to the Segel-Jackson model. Near the Bautin bifurcation we theoretically proved that solutions with small (large) initial values are convergent to stable periodic oscillations (diverge to infinity).

\section*{Acknowledgments} The authors wish to express their gratitude to the editors
and the reviewers for the helpful comments. This research is supported by  National Natural Science Foundation of China (11701120,11771109).

\section*{Appendix}
In this Appendix, we will give  detailed formulae to calculating $g_{ij}$'s, then the first and second Lyapunov coefficients are determined.

In fact, by comparing the same order of terms, we have
$$g_{20}=2\bar{M}\tau^\ast\left[\bar{q_2}(k-a q_1)+(-q_2^2+ q_1\textrm{e}^{-2\textrm{i}\omega^\ast\tau^\ast})\right]\int_0^{l\pi}b_{n_0}^3d x,$$
$$g_{11}=\bar{M}\tau^\ast\{\bar{q_2}[2 k-a(\bar{q_1} q_1)]-2 q_1\bar{q_1}+\bar{q_1}+q_1\}\int_0^{l\pi}b_{n_0}^3d x,$$
$$g_{02}=2\bar{M}\tau^\ast\left[\bar{q_2}(k-a \bar{q_1})+(-\bar{q_2}^2+ \bar{q_1}\textrm{e}^{2\textrm{i}\omega^\ast\tau^\ast})\right]\int_0^{l\pi}b_{n_0}^3d x,$$
$$g_{30}=3\bar{M}\tau^\ast\int_0^{l\pi}\left\{\bar{q_2}(k+1)W_{20}^{(1)}(0)-q_1\left[W_{20}^{(1)}(0)+W_{20}^{(2)}(0)\right]+2W_{20}^{(1)}(-1)\right\}b_{n_0}^2d x,$$
$$\begin{array}{ll}g_{21}=&\bar{M}\tau^\ast\bar{q_2}\int_0^{l\pi}\left\{2k\left[W_{11}^{(1)}(0)+W_{20}^{(1)}(0)\right]\right.\\& \left.-a\left[2W_{11}^{(2)}(0)+W_{20}^{(2)}(0)+\bar{q_1}W_{20}^{(1)}(0)+2q_1W_{11}^{(1)}(0)\right]\right\}b_{n_0}^2d x\\&+\frac12\bar M\tau^\ast\int_0^{l\pi}\left\{-\left[8q_1W_{11}^{(2)}(0)+2\bar q_1W_{20}^{(2)}(0)+W_{20}^{(2)}(0)W_{02}^{(2)}(0)\right]+\left[4\textrm{e}^{-\textrm{i}\omega^\ast\tau^\ast}W_{11}^{(2)}(-1)\right.\right.\\&\left.\left.+2\textrm{e}^{\textrm{i}\omega^\ast\tau^\ast}W_{20}^{(2)}(-1)+2\textrm{e}^{\textrm{i}\omega^\ast\tau^\ast}\bar q_1W_{20}^{(1)}(-1)+4\textrm{e}^{-\textrm{i}\omega^\ast\tau^\ast}q_1W_{11}^{(1)}(-1)\right]\right\}b_{n_0}^2d x,\end{array}$$
$$\begin{array}{ll}g_{12}=&\bar{M}\tau^\ast\bar{q_2}\int_0^{l\pi}\left\{k\left[3W_{02}^{(1)}(0)+4W_{11}^{(1)}(0)\right]\right.\\& \left.-a\left[2W_{11}^{(2)}(0)+W_{02}^{(2)}(0)+{q_1}W_{02}^{(1)}(0)+2\bar q_1W_{11}^{(1)}(0)\right]\right\}b_{n_0}^2d x\\&+\bar M\tau^\ast\int_0^{l\pi}\left\{-\left[(q_1+\bar q_1W_{20}^{(2)}(0)+2(\bar q_1+1)W_{11}^{(2)}(0)+q_1W_{02}^{(2)}(0)\right]+\left[\textrm{e}^{-\textrm{i}\omega^\ast\tau^\ast}W_{02}^{(2)}(-1)\right.\right.\\&\left.\left.+2\textrm{e}^{\textrm{i}\omega^\ast\tau^\ast}W_{11}^{(2)}(-1)+2\bar q_1\textrm{e}^{\textrm{i}\omega^\ast\tau^\ast}W_{11}^{(1)}(-1)+\textrm{e}^{-\textrm{i}\omega^\ast\tau^\ast}q_1W_{02}^{(1)}(-1)\right]\right\}b_{n_0}^2d x,\end{array}$$
$$\begin{array}{ll}g_{03}=&3\bar{M}\tau^\ast\int_0^{l\pi}\left\{2\bar q_2 kW_{20}^{(1)}(0)-a\bar q_2\left[W_{02}^{(2)}(0)+\bar q_1W_{02}^{(1)}(0)\right]-2\bar q_1W_{02}^{(2)}(0)\right.\\ &+\left. \textrm{e}^{-\textrm{i}\omega^\ast\tau^\ast}\left[W_{02}^{(2)}(-1)+\bar q_1W_{02}^{(1)}(-1)\right]\right\}b_{n_0}^2d x,\end{array}$$
$$\begin{array}{ll}g_{40}=&2\bar{M}\tau^\ast\bar q_2 \left\{k\left[4\int_0^{l\pi}W_{30}^{(1)}(0)b_{n_0}^2dx+3\int_0^{l\pi}\left(W_{20}^{(1)}(0)\right)^2b_{n_0}dx\right]\right.\\&\left.-2a\int_0^{l\pi}\left[W_{30}^{(2)}(0)+q_1W_{30}^{(1)}(0)\right]b_{n_0}^2dx-3a\int_0^{l\pi} W_{20}^{(1)}(0)W_{20}^{(2)}(0) b_{n_0}dx\right\}\\
&+2\bar M\tau^\ast\left\{-\left[4q_1\int_0^{l\pi}W_{30}^{(2)}(0)b_{n_0}^2dx+3\int_0^{l\pi}\left[2\left(W_{20}^{(2)}(0)\right)^2+W_{20}^{(1)}(-1)W_{20}^{(2)}(-1)\right]b_{n_0}dx\right]\right.\\
&\left.+2\textrm{e}^{-\textrm{i}\omega^\ast\tau^\ast}\int_0^{l\pi}\left[W_{30}^{(2)}(-1)+q_1W_{30}^{(1)}(-1)\right]b_{n_0}^2dx\right\},\end{array}$$
$$\begin{array}{ll}g_{31}=&\bar{M}\tau^\ast\bar q_2 \left\{k\int_0^{l\pi}\left[9W_{21}^{(1)}(0)+W_{30}^{(1)}(0) \right]b_{n_0}^2dx+\int_0^{l\pi}\left[k-3aW_{20}^{(2)}(0)\right]W_{20}^{(1)}(0)W_{11}^{(1)}(0)b_{n_0}dx\right.\\&\left.-a\int_0^{l\pi}\left[3W_{21}^{(2)}(0)+W_{30}^{(2)}(0)+\bar q_1W_{30}^{(1)}(0)+3q_1W_{21}^{(1)}(0) \right]b_{n_0}^2dx\right\}\\
&+\bar M\tau^\ast\left\{-\int_0^{l\pi}\left[3q_1W_{21}^{(2)}(0)+2\bar q_1W_{30}^{(2)}(0)\right]b_{n_0}^2dx\right.\\&\left.-\int_0^{l\pi}\left[W_{20}^{(2)}(0)W_{11}^{(2)}(0)-W_{20}^{(1)}(-1)W_{11}^{(2)}(-1)\right]b_{n_0}dx\right.\\
&\left.+3\textrm{e}^{-\textrm{i}\omega^\ast\tau^\ast}\int_0^{l\pi}\left[W_{21}^{(2)}(-1)+q_1W_{21}^{(1)}(-1)\right]b_{n_0}^2dx\right.\\&\left.+\textrm{e}^{\textrm{i}\omega^\ast\tau^\ast}\int_0^{l\pi}\left[W_{30}^{(2)}(-1)+\bar q_1W_{31}^{(1)}(-1)\right]b_{n_0}^2dx\right\},\end{array}$$
$$\begin{array}{ll}g_{22}=&\bar{M}\tau^\ast\bar q_2 \left\{4k\int_0^{l\pi}\left[W_{12}^{(1)}(0)+W_{21}^{(1)}(0) \right]b_{n_0}^2dx+k\int_0^{l\pi}\left[2W_{20}^{(1)}(0)W_{02}^{(1)}(0)+\left(W_{11}^{(1)}(0)\right)^2\right]b_{n_0}dx\right.\\&\left.-2a\int_0^{l\pi}\left[W_{12}^{(2)}(0)+W_{20}^{(2)}(0)+\bar q_1W_{21}^{(1)}(0)+q_1W_{12}^{(1)}(0) \right]b_{n_0}^2dx\right.\\
&\left.-a\int_0^{l\pi}\left[W_{20}^{(1)}(0)W_{02}^{(2)}(0)+4W_{11}^{(1)}(0)W_{11}^{(2)}(0)+W_{02}^{(1)}(0)W_{20}^{(2)}(0)\right]b_{n_0}dx\right\}\\
&+\bar M\tau^\ast\left\{-4\int_0^{l\pi}\left[q_1W_{12}^{(2)}(0)+\bar q_1W_{21}^{(2)}(0)\right]b_{n_0}^2dx+4\left(\textrm{e}^{-\textrm{i}\omega^\ast\tau^\ast}+\textrm{e}^{\textrm{i}\omega^\ast\tau^\ast}   \right)\int_0^{l\pi}W_{21}^{(2)}(-1)b_{n_0}^2dx\right.\\
&\left.-\int_0^{l\pi}\left[W_{20}^{(2)}(0)W_{02}^{(2)}(0)+\left(W_{11}^{(2)}(0)\right)^2-W_{11}^{(1)}(-1)W_{11}^{(2)}(-1)\right]b_{n_0}dx\right.\\
&\left.+2\int_0^{l\pi}\left[\bar q_1 \textrm{e}^{\textrm{i}\omega^\ast\tau^\ast}W_{21}^{(1)}(-1)+q_1\textrm{e}^{-\textrm{i}\omega^\ast\tau^\ast}W_{12}^{(1)}(-1)\right]b_{n_0}^2dx\right.\\&\left.+\int_0^{l\pi}W_{22}^{(2)}(-1)\left[W_{20}^{(1)}(-1)+W_{02}^{(1)}(-1)\right]b_{n_0}dx\right\},\end{array}$$
$$\begin{array}{ll}g_{13}=&\bar{M}\tau^\ast\bar q_2 \left\{2k\int_0^{l\pi}\left[W_{03}^{(1)}(0)+3W_{12}^{(1)}(0) \right]b_{n_0}^2dx+k\int_0^{l\pi}W_{02}^{(1)}(0)W_{11}^{(1)}(0)b_{n_0}dx\right.\\&\left.-a\int_0^{l\pi}\left[3W_{12}^{(2)}(0)+W_{03}^{(2)}(0)+ q_1W_{03}^{(1)}(0)+3\bar q_1W_{12}^{(1)}(0) \right]b_{n_0}^2dx\right.\\&\left.-3a\int_0^{l\pi}\left[W_{11}^{(1)}(0)W_{02}^{2)}(0)+W_{02}^{(1)}(0)W_{11}^{(2)}(0) \right]b_{n_0}dx\right\}\\
&+\bar M\tau^\ast\left\{-2\int_0^{l\pi}\left[q_1W_{03}^{2)}(0)+3\bar q_1W_{12}^{2)}(0)\right]b_{n_0}^2dx\right.\\&\left.+\textrm{e}^{-\textrm{i}\omega^\ast\tau^\ast}\int_0^{l\pi}\left[W_{03}^{(2)}(-1)+q_1W_{03}^{(1)}(-1)\right]b_{n_0}^2dx\right.\\&\left.+3\textrm{e}^{\textrm{i}\omega^\ast\tau^\ast}\int_0^{l\pi}\left[W_{12}^{(2)}(-1)+\bar q_1W_{12}^{(1)}(-1)\right]b_{n_0}^2dx\right.\\&\left.+3\int_0^{l\pi}W_{02}^{(1)}(-1)W_{11}^{(2)}(-1)b_{n_0}dx\right\},\end{array}$$
$$\begin{array}{ll}g_{04}=&2\bar{M}\tau^\ast\bar q_2 \left\{\int_0^{l\pi}\left[(14k-2a\bar q_1)W_{03}^{(1)}(0)-2aW_{03}^{(2)}(0)\right]b_{n_0}^2dx\right.\\&\left.+3(k-a)\int_0^{l\pi}W_{02}^{(1)}(0)\left[W_{02}^{(1)}(0)+W_{02}^{(2)}(0)\right]b_{n_0}dx\right\}\\
&+2\bar M\tau^\ast\left\{-4\bar q_1\int_0^{l\pi}W_{03}^{(2)}(0)b_{n_0}^2dx+2\textrm{e}^{\textrm{i}\omega^\ast\tau^\ast}\int_0^{l\pi}\left[W_{03}^{(2)}(-1)+\bar q_1W_{03}^{(1)}(-1)\right]b_{n_0}^2dx\right.\\
&\left.-3\int_0^{l\pi}\left[\left(W_{02}^{(2)}(0)\right)^2-W_{02}^{(1)}(-1)W_{02}^{(2)}(-1)\right]b_{n_0}dx\right\}\end{array}$$
and$$\begin{array}{ll}g_{32}=&\bar{M}\tau^\ast\bar q_2 \left\{2k\int_0^{l\pi}\left[3W_{22}^{(1)}(0)+2W_{31}^{(1)}(0)\right]b_{n_0}^2dx\right.\\&\left.+k\int_0^{l\pi}\left[9W_{20}^{(1)}(0)W_{12}^{(1)}(0)+12W_{11}^{(1)}(0)W_{21}^{(1)}(0)+2W_{02}^{(1)}(0)W_{30}^{(1)}(0)\right]b_{n_0}dx\right.\\
&\left.-a\int_0^{l\pi}\left[3W_{22}^{(2)}(0)+2W_{31}^{(2)}(0)+2\bar q_1 W_{31}^{(1)}(0)+3q_1W_{21}^{(1)}(0)\right]b_{n_0}^2dx\right.\\
&\left.-a\int_0^{l\pi}\left[3W_{20}^{(1)}(0)W_{12}^{(2)}(0)+6W_{11}^{(1)}(0)W_{21}^{(2)}(0)+W_{02}^{(1)}(0)W_{30}^{(2)}(0)\right.\right.\\
&\left.\left.+W_{30}^{(1)}(0)W_{02}^{(2)}(0)+6W_{21}^{(1)}(0)W_{11}^{(2)}(0)+3W_{12}^{(1)}(0)W_{20}^{(2)}(0)\right]b_{n_0}dx\right\}\\
&+2\bar M\tau^\ast\left\{-2\int_0^{l\pi}\left[3q_1W_{22}^{(2)}(0)+2\bar q_1W_{31}^{(2)}(0)\right]b_{n_0}^2dx\right.\\&\left.-\int_0^{l\pi}\left[12W_{11}^{(2)}(0)W_{21}^{(2)}(0)+2W_{02}^{(2)}(0)W_{30}^{(2)}(0)+3W_{12}^{(2)}(0)W_{20}^{(2)}(0)\right]b_{n_0}dx\right.\\
&\left.+3\textrm{e}^{-\textrm{i}\omega^\ast\tau^\ast}\int_0^{l\pi}\left[W_{22}^{(2)}(-1)+q_1W_{22}^{(1)}(-1)\right]b_{n_0}^2dx\right.\\
&\left.+2\textrm{e}^{\textrm{i}\omega^\ast\tau^\ast}\int_0^{l\pi}\left[W_{31}^{(2)}(-1)+\bar q_1W_{31}^{(1)}(-1)\right]b_{n_0}^2dx\right.\\
&\left.+\int_0^{l\pi}\left[3W_{20}^{(1)}(-1)W_{12}^{(2)}(-1)+W_{02}^{(1)}(-1)W_{30}^{(2)}(-1)\right.\right.\\
&\left.\left.+W_{30}^{(1)}(-1)W_{02}^{(2)}(-1)+6W_{21}^{(1)}(-1)W_{11}^{(2)}(-1)\right]b_{n_0}dx\right\}.\end{array}$$

In the current step, we still need to calculate $W_{ij}(\theta)$. Express $H(z,\bar{z},\theta)$ by the Taylor series, we have
$$H(z,\bar{z},\cdot)=H_{20}\frac{z^2}{2}+H_{11} z\bar{z}+H_{02}\frac{\bar{z}^2}{2}+H_{30}\frac{z^3}{6}+H_{21}\frac{z^2\bar{z}}{2}+H_{12}\frac{z\bar{z}^2}{2}+H_{03}\frac{\bar{z}^3}{6}+H_{31}\frac{z^3\bar{z}}{6}+H_{22}\frac{z^2\bar{z}^2}{4}+\cdots.$$

Obviously
$$H_{i j}(\theta)=\left\{\begin{array}{ll}
-g_{i j}q(\theta)b_{n_0}-\bar{g}_{j i}\bar{q}(0)b_{n_0} ,& \theta\in[-r,0),\\
 -g_{i j}q(0)b_{n_0}-\bar{g}_{j i}\bar{q}(0)b_{n_0}+\frac{\partial^{(i+j)}\tilde F}{\partial^i z\partial^j \bar z},& \theta=0.\end{array}
 \right.$$

Notice the definition of $A_0$ and Eq.(\ref{Waij}), then we have \begin{equation}\label{Wij}
 \begin{array}{l}
 (A_0-2{\rm i}\omega^\ast\tau^\ast I)W_{20}(\theta)=-H_{20}(\theta),~A_0W_{11}(\theta)=-H_{11}(\theta),\\(A_0+2{\rm i}\omega^\ast\tau^\ast I)W_{02}(\theta)=-H_{02}(\theta),\\
 (A_0-3{\rm i}\omega^\ast\tau^\ast I)W_{30}(\theta)-3g_{20}W_{20}(\theta)-3\bar g_{02}W_{11}(\theta)=-H_{30}(\theta),\\
 (A_0-{\rm i}\omega^\ast\tau^\ast I)W_{21}(\theta)-(g_{20}+2\bar g_{11})W_{11}(\theta)-2 g_{11}W_{20}(\theta)-\bar g_{02}W_{02}(\theta)=-H_{21}(\theta),\\
  (A_0+{\rm i}\omega^\ast\tau^\ast I)W_{12}(\theta)-(\bar g_{20}+2 g_{11})W_{11}(\theta)-2\bar g_{11}W_{02}(\theta)- g_{02}W_{20}(\theta)=-H_{12}(\theta),\\
   (A_0+3{\rm i}\omega^\ast\tau^\ast I)W_{30}(\theta)-3g_{02}W_{11}(\theta)-3\bar g_{20}W_{02}(\theta)=-H_{03}(\theta),\\
   (A_0-2{\rm i}\omega^\ast\tau^\ast I)W_{31}(\theta)-3(g_{20}+\bar g_{11})W_{21}(\theta)-(g_{30}+3\bar g_{12})W_{11}(\theta)-3  g_{11}W_{30}(\theta)\\~~~~-3g_{21}W_{20}(\theta)-\bar g_{03}W_{02}(\theta)-3\bar g_{02}W_{12}(\theta)=-H_{31}(\theta),\\
   A_0W_{22}(\theta)-g_{02}W_{30}(\theta)-(g_{20}+4\bar g_{11})W_{12}(\theta)-2(g_{21}+\bar g_{21})W_{11}(\theta)-2g_{12}W_{20}(\theta)
   \\~~~~-(\bar g_{20}+4g_{11})W_{21}(\theta)-2\bar g_{12}W_{02}(\theta)-\bar g_{02} W_{03}(\theta)=-H_{22}(\theta).
 \end{array}
 \end{equation}

Similar with the results given in \cite{Xu1}, we can solve these functions by
$$W_{20}(\theta)=\frac{{\rm i} g_{20}}{\omega^\ast\tau^\ast}q(0)b_{n_0}{\rm e}^{{\rm i}\omega^\ast\tau^\ast\theta}+\frac{{\rm i} \bar g_{20}}{3\omega^\ast\tau^\ast}\bar q(0)b_{n_0}{\rm e}^{-{\rm i}\omega^\ast\tau^\ast\theta}+T_1{\rm e}^{2{\rm i}\omega^\ast\tau^\ast\theta}, $$
$$W_{11}(\theta)=-\frac{{\rm i} g_{11}}{\omega^\ast\tau^\ast}q(0)b_{n_0}{\rm e}^{{\rm i}\omega^\ast\tau^\ast\theta}+\frac{{\rm i} \bar g_{11}}{\omega^\ast\tau^\ast}\bar q(0)b_{n_0}{\rm e}^{-{\rm i}\omega^\ast\tau^\ast\theta}+T_2, $$
$$W_{02}(\theta)=-\frac{{\rm i} \bar g_{20}}{\omega^\ast\tau^\ast}\bar q(0)b_{n_0}{\rm e}^{-{\rm i}\omega^\ast\tau^\ast\theta}-\frac{{\rm i}  g_{02}}{3\omega^\ast\tau^\ast} q(0)b_{n_0}{\rm e}^{{\rm i}\omega^\ast\tau^\ast\theta}+T_3{\rm e}^{-2{\rm i}\omega^\ast\tau^\ast\theta}, $$
$$W_{30}(\theta)=\frac{{\rm i} S_{30}^{(1)}}{2\omega^\ast\tau^\ast}q(0)b_{n_0}{\rm e}^{{\rm i}\omega^\ast\tau^\ast\theta}+\frac{{\rm i} S_{30}^{(2)}}{4\omega^\ast\tau^\ast}\bar q(0)b_{n_0}{\rm e}^{-{\rm i}\omega^\ast\tau^\ast\theta}+\frac{{\rm i} S_{30}^{(3)}}{\omega^\ast\tau^\ast}\bar q(0)b_{n_0}{\rm e}^{2{\rm i}\omega^\ast\tau^\ast\theta}+\frac{{\rm i} S_{30}^{(4)}}{3\omega^\ast\tau^\ast}+T_4{\rm e}^{3{\rm i}\omega^\ast\tau^\ast\theta}, $$
$$W_{21}(\theta)=S_{21}^{(1)} q(0)b_{n_0}\theta{\rm e}^{{\rm i}\omega^\ast\tau^\ast\theta}+\frac{{\rm i} S_{21}^{(2)}}{2\omega^\ast\tau^\ast}\bar q(0)b_{n_0}{\rm e}^{-{\rm i}\omega^\ast\tau^\ast\theta}-\frac{{\rm i} S_{21}^{(3)}}{\omega^\ast\tau^\ast}{\rm e}^{2{\rm i}\omega^\ast\tau^\ast\theta}+\frac{{\rm i} S_{21}^{(4)}}{3\omega^\ast\tau^\ast}{\rm e}^{-2{\rm i}\omega^\ast\tau^\ast\theta}+\frac{{\rm i} S_{21}^{(5)}}{\omega^\ast\tau^\ast}+T_5{\rm e}^{{\rm i}\omega^\ast\tau^\ast\theta}, $$
$$W_{12}(\theta)=-\frac{{\rm i} S_{12}^{(1)}}{2\omega^\ast\tau^\ast} q(0)b_{n_0}{\rm e}^{{\rm i}\omega^\ast\tau^\ast\theta}+S_{12}^{(2)} \bar q(0)b_{n_0}\theta{\rm e}^{-{\rm i}\omega^\ast\tau^\ast\theta}-\frac{{\rm i} S_{12}^{(3)}}{3\omega^\ast\tau^\ast}{\rm e}^{2{\rm i}\omega^\ast\tau^\ast\theta}+\frac{{\rm i} S_{12}^{(4)}}{\omega^\ast\tau^\ast}{\rm e}^{-2{\rm i}\omega^\ast\tau^\ast\theta}-\frac{{\rm i} S_{12}^{(5)}}{\omega^\ast\tau^\ast}+T_6{\rm e}^{-{\rm i}\omega^\ast\tau^\ast\theta}, $$
$$W_{03}(\theta)=-\frac{{\rm i} S_{03}^{(1)}}{4\omega^\ast\tau^\ast} q(0)b_{n_0}{\rm e}^{{\rm i}\omega^\ast\tau^\ast\theta}-\frac{{\rm i} S_{03}^{(2)}}{2\omega^\ast\tau^\ast}\bar q(0)b_{n_0}{\rm e}^{-{\rm i}\omega^\ast\tau^\ast\theta}-\frac{{\rm i} S_{03}^{(3)}}{\omega^\ast\tau^\ast}{\rm e}^{-2{\rm i}\omega^\ast\tau^\ast\theta}-\frac{{\rm i} S_{03}^{(4)}}{3\omega^\ast\tau^\ast}+T_7{\rm e}^{-3{\rm i}\omega^\ast\tau^\ast\theta},$$
$$\begin{array}{l}W_{31}(\theta)=\frac{{\rm i} S_{31}^{(1)}}{\omega^\ast\tau^\ast} {\rm e}^{{\rm i}\omega^\ast\tau^\ast\theta}+\frac{S_{31}^{(11)}}{(\omega^\ast\tau^\ast)^2}{\rm e}^{{\rm i}\omega^\ast\tau^\ast\theta}+\frac{{\rm i} S_{31}^{(11)}}{\omega^\ast\tau^\ast}\theta{\rm e}^{{\rm i}\omega^\ast\tau^\ast\theta}+\frac{{\rm i} S_{31}^{(2)}}{3\omega^\ast\tau^\ast}{\rm e}^{-{\rm i}\omega^\ast\tau^\ast\theta}+\frac{S_{31}^{(21)}}{9(\omega^\ast\tau^\ast)^2}{\rm e}^{-{\rm i}\omega^\ast\tau^\ast\theta}\\+\frac{{\rm i} S_{31}^{(21)}}{3\omega^\ast\tau^\ast}\theta {\rm e}^{-{\rm i}\omega^\ast\tau^\ast\theta}+ S_{31}^{(3)}\theta {\rm e}^{2{\rm i}\omega^\ast\tau^\ast\theta}+\frac{{\rm i} S_{31}^{(4)}}{4\omega^\ast\tau^\ast}{\rm e}^{-2{\rm i}\omega^\ast\tau^\ast\theta}-\frac{{\rm i} S_{31}^{(5)}}{\omega^\ast\tau^\ast}{\rm e}^{3{\rm i}\omega^\ast\tau^\ast\theta}+\frac{{\rm i} S_{31}^{(6)}}{2\omega^\ast\tau^\ast}+T_8{\rm e}^{2{\rm i}\omega^\ast\tau^\ast\theta} \end{array}$$and
$$\begin{array}{l}W_{22}(\theta)=-\frac{{\rm i} S_{22}^{(1)}}{\omega^\ast\tau^\ast} {\rm e}^{{\rm i}\omega^\ast\tau^\ast\theta}-\frac{S_{22}^{(11)}}{\omega^\ast\tau^\ast}\theta{\rm e}^{{\rm i}\omega^\ast\tau^\ast\theta}+\frac{{\rm i} S_{22}^{(11)}}{(\omega^\ast\tau^\ast)^2}{\rm e}^{{\rm i}\omega^\ast\tau^\ast\theta}+\frac{{\rm i} S_{22}^{(2)}}{\omega^\ast\tau^\ast}{\rm e}^{-{\rm i}\omega^\ast\tau^\ast\theta}+\frac{S_{22}^{(21)}}{\omega^\ast\tau^\ast}\theta{\rm e}^{-{\rm i}\omega^\ast\tau^\ast\theta}\\+\frac{{\rm i} S_{22}^{(21)}}{(\omega^\ast\tau^\ast)^2} {\rm e}^{-{\rm i}\omega^\ast\tau^\ast\theta}-\frac{ {\rm i} S_{22}^{(3)}}{2\omega^\ast\tau^\ast}{\rm e}^{2{\rm i}\omega^\ast\tau^\ast\theta}+\frac{{\rm i} S_{22}^{(4)}}{2\omega^\ast\tau^\ast}{\rm e}^{-2{\rm i}\omega^\ast\tau^\ast\theta}-\frac{{\rm i} S_{22}^{(5)}}{3\omega^\ast\tau^\ast}{\rm e}^{3{\rm i}\omega^\ast\tau^\ast\theta}+\frac{{\rm i} S_{22}^{(6)}}{3\omega^\ast\tau^\ast}{\rm e}^{-3{\rm i}\omega^\ast\tau^\ast\theta}+S_{22}^{(7)}\theta+T_9. \end{array}$$

The coefficients with form $S_{i j}^{(k)}$'s are given by
$$S_{30}^{(1)}=3\frac{{\rm i} g_{20}^2}{\omega^\ast\tau^\ast}+3\frac{{\rm i} g_{20} \bar g_{02}}{\omega^\ast\tau^\ast}+g_{30},~S_{30}^{(2)}=\frac{{\rm i}\bar  g_{02} g_{20}}{\omega^\ast\tau^\ast}+3\frac{{\rm i} \bar g_{11} \bar g_{02}}{\omega^\ast\tau^\ast}+\bar g_{03},~S_{30}^{(3)}=3g_{20}T_1,~S{30}^{(4)}=3\bar g_{02}T_2,$$
$$S_{21}^{(1)}=-(g_{20}+2\bar g_{11})\frac{{\rm i} g_{11}}{\tau^\ast\omega^\ast}+2g_{11}\frac{{\rm i} g_{20}}{\tau^\ast\omega^\ast}-\frac{{\rm i} g_{02}\bar g_{02}}{3\tau^\ast\omega^\ast}+g_{21},$$
$$S_{21}^{(2)}=(g_{20}+2\bar g_{11})\frac{{\rm i} \bar g_{11}}{\tau^\ast\omega^\ast}+2g_{11}\frac{{\rm i} \bar g_{02}}{3\tau^\ast\omega^\ast}-\frac{{\rm i} \bar g_{20}\bar g_{02}}{\tau^\ast\omega^\ast}+\bar g_{12},$$
$$S_{21}^{(3)}=2g_{11}T_1,S_{21}^{(4)}=\bar g_{02}T_3,S_{21}^{(5)}=(g_{02}+2\bar g_{11})T_2,$$
$$S_{12}^{(1)}=\frac{{\rm i} g_{20}g_{02}}{\omega^\ast\tau^\ast}-(2g_{11}+\bar g_{20})\frac{{\rm i} g_{11} }{\omega^\ast\tau^\ast}-2\bar g_{11}\frac{{\rm i} g_{02} }{3\omega^\ast\tau^\ast}+g_{12},$$
$$S_{12}^{(2)}=\frac{{\rm i} \bar g_{02}g_{02}}{3\omega^\ast\tau^\ast}+(2g_{11}+\bar g_{20})\frac{{\rm i} g_{11} }{\omega^\ast\tau^\ast}-2\frac{{\rm i} \bar g_{20}\bar g_{11} }{\omega^\ast\tau^\ast}+\bar g_{21},$$
$$S_{12}^{(3)}=g_{02}T_1,S_{12}^{(4)}=2\bar g_{11}T_3,S_{12}^{(5)}=(2g_{11}+2\bar g_{20})T_2,$$
$$S_{03}^{(1)}=-3\frac{{\rm i} g_{11}g_{02}}{\omega^\ast\tau^\ast}-\frac{{\rm i} g_{02} \bar g_{20}}{\omega^\ast\tau^\ast}+g_{03},~S_{03}^{(2)}=3\frac{{\rm i}\bar  g_{11} g_{02}}{\omega^\ast\tau^\ast}-3\frac{{\rm i} \bar g_{20}^2 }{\omega^\ast\tau^\ast}+\bar g_{30},~S_{03}^{(3)}=3\bar g_{20}T_3,~S_{03}^{(4)}=3g_{02}T_2,$$
$$S_{31}^{(1)}=\left[3g_{21}\frac{{\rm i} g_{20}}{\omega^\ast\tau^\ast}-(g_{30}+3\bar g_{12})\frac{{\rm i} g_{11}}{\omega^\ast\tau^\ast}-\frac{{\rm i} g_{02}\bar g_{03}}{3\omega^\ast\tau^\ast}+3g_{11}\frac{{\rm i} S_{30}^{(1)}}{2\omega^\ast\tau^\ast}-3\bar g_{02}\frac{{\rm i} S_{12}^{(1)}}{2\omega^\ast\tau^\ast}+g_{31}\right]q(0)b_{n_0}+3(g_{20}+\bar g_{11})T_5,$$
$$S_{31}^{(11)}=3(g_{20}+\bar g_{11})S_{21}^{(1)}q(0)b_{n_0},$$
$$S_{31}^{(2)}=\left[g_{21}\frac{{\rm i} \bar g_{02}}{\omega^\ast\tau^\ast}+(g_{30}+3\bar g_{12})\frac{{\rm i}\bar g_{11}}{\omega^\ast\tau^\ast}-\frac{{\rm i} \bar g_{20}\bar g_{03}}{\omega^\ast\tau^\ast}+3g_{11}\frac{{\rm i} S_{30}^{(2)}}{4\omega^\ast\tau^\ast}+3(\bar g_{11}+g_{20})\frac{{\rm i} S_{21}^{(2)}}{2\omega^\ast\tau^\ast}+\bar g_{13}\right]q(0)b_{n_0}-3g_{02}T_6,$$
$$S_{31}^{(21)}=3\bar g_{02}S_{12}^{(2)}\bar q(0)b_{n_0},$$
$$S_{31}^{(3)}=3g_{21}T_1+3 g_{11}\frac{{\rm i} S_{30}^{(3)}}{\omega^\ast\tau^\ast}-3(g_{20}+\bar g_{11})\frac{{\rm i} S_{21}^{(3)}}{\omega^\ast\tau^\ast}-3\bar g_{02}\frac{{\rm i} S_{12}^{(3)}}{3\omega^\ast\tau^\ast},$$
$$S_{31}^{(4)}=\bar g_{03}T_3+ (\bar g_{11}+g_{20})\frac{{\rm i} S_{21}^{(4)}}{\omega^\ast\tau^\ast}+3\bar g_{02}\frac{{\rm i} S_{12}^{(4)}}{\omega^\ast\tau^\ast},$$
$$S_{31}^{(5)}=3g_{11}T_4,$$
$$S_{31}^{(6)}=(g_{30}+3\bar g_{12})T_2+  g_{11} \frac{{\rm i} S_{30}^{(4)}}{\omega^\ast\tau^\ast}+3(\bar g_{11}+g_{20})\frac{{\rm i} S_{21}^{(5)}}{\omega^\ast\tau^\ast}-3\bar g_{02}\frac{{\rm i} S_{12}^{(5)}}{\omega^\ast\tau^\ast},$$
$$\begin{array}{l}S_{22}^{(1)}=\left[g_{02}\frac{{\rm i} S_{30}^{(1)}}{2\omega^\ast\tau^\ast}-(g_{20}+4\bar g_{11})\frac{{\rm i} S_{12}^{(1)}}{2\omega^\ast\tau^\ast}-2(g_{21}+\bar g_{21})\frac{{\rm i} g_{11}}{\omega^\ast\tau^\ast}+2g_{12} \frac{{\rm i} g_{20}}{\omega^\ast\tau^\ast}-2\bar g_{12} \frac{{\rm i} g_{02}}{3\omega^\ast\tau^\ast}-\bar g_{02}\frac{{\rm i} S_{03}^{(1)}}{4\omega^\ast\tau^\ast}+g_{22}\right]\\
~~~~\times q(0)b_{n_0}+(4g_{11}+\bar g_{20})T_5,\end{array}$$
$$S_{22}^{(11)}=(4g_{11}+\bar g_{20})S_{21}^{(1)}q(0)b_{n_0},$$
$$\begin{array}{l}S_{22}^{(2)}=\left[g_{02}\frac{{\rm i} S_{30}^{(2)}}{4\omega^\ast\tau^\ast}+2(g_{21}+\bar g_{21})\frac{{\rm i} \bar g_{11}}{\omega^\ast\tau^\ast}+2g_{12}\frac{{\rm i} \bar g_{02}}{3\omega^\ast\tau^\ast}+4(g_{11}+\bar g_{20}) \frac{{\rm i} S_{21}^{(2)}}{2\omega^\ast\tau^\ast}-2\bar g_{12} \frac{{\rm i} \bar g_{20}}{\omega^\ast\tau^\ast}-\bar g_{02}\frac{{\rm i} S_{03}^{(2)}}{2\omega^\ast\tau^\ast}+\bar g_{22}\right]\\
~~~~\times \bar q(0)b_{n_0}+(4\bar g_{11}+ g_{20})T_6,\end{array}$$
$$S_{22}^{(21)}=(4\bar g_{11}+ g_{20})S_{12}^{(2)}\bar q(0)b_{n_0},$$
$$S_{22}^{(3)}=[g_{02}\frac{{\rm i} S_{30}^{(3)}}{\omega^\ast\tau^\ast}-(g_{20}+4\bar g_{11})\frac{{\rm i} S_{12}^{(3)}}{3\omega^\ast\tau^\ast}+2g_{12}T_1-(\bar g_{20}+4 g_{11})\frac{{\rm i} S_{21}^{(3)}}{\omega^\ast\tau^\ast},$$
$$S_{22}^{(4)}=-\bar g_{02}\frac{{\rm i} S_{03}^{(3)}}{\omega^\ast\tau^\ast}+(g_{20}+4\bar g_{11})\frac{{\rm i} S_{12}^{(4)}}{\omega^\ast\tau^\ast}+2\bar g_{12}T_3+(\bar g_{20}+4 g_{11})\frac{{\rm i} S_{21}^{(4)}}{3\omega^\ast\tau^\ast},$$
$$S_{22}^{(5)}=g_{02}T_4, ~~~S_{22}^{(6)}=\bar g_{02}$$ and
$$S_{22}^{(7)}=g_{02}\frac{{\rm i} S_{30}^{(4)}}{3\omega^\ast\tau^\ast}-(g_{20}+4\bar g_{11})\frac{{\rm i} S_{21}^{(5)}}{\omega^\ast\tau^\ast}+2(g_{21}+\bar g_{21})T_2 +4(g_{11}+\bar g_{20}) \frac{{\rm i} S_{21}^{(5)}}{\omega^\ast\tau^\ast}-\bar g_{02} \frac{{\rm i} S_{03}^{(4)}}{3\omega^\ast\tau^\ast}.$$

Now the only unknowns are $T_j's$, they can be calculated by expanding $T_j$ as $T_j=\sum_{n=0}^\infty T_j^nb_n$, then from (\ref{Wij}), we have

$$\left(2{\rm i}\omega^\ast\tau^\ast I-\int_{-1}^0d\eta_n(0,\theta){\rm e}^{2{\rm i}\omega^\ast\tau^\ast\theta}\right)T_1^n=<\tilde F''_{zz},\beta_n>,$$
$$\int_{-1}^0d\eta_n(0,\theta)T_2^n=-<\tilde F''_{z\bar z},\beta_n>,$$
$$\left(-2{\rm i}\omega^\ast\tau^\ast I-\int_{-1}^0d\eta_n(0,\theta){\rm e}^{-2{\rm i}\omega^\ast\tau^\ast\theta}\right)T_3^n=<\tilde F''_{\bar z\bar z},\beta_n>,$$
$$\left(3{\rm i}\omega^\ast\tau^\ast I-\int_{-1}^0d\eta_n(0,\theta){\rm e}^{3{\rm i}\omega^\ast\tau^\ast\theta}\right)T_4^n=<\tilde F'''_{zzz}-\frac{3{\rm i} g_{20}}{\omega^\ast\tau^\ast}\tilde F''_{zz}-\frac{{\rm i} \bar g_{02}}{\omega^\ast\tau^\ast}\tilde F''_{z\bar z},\beta_n>,$$
$$\begin{array}{l}\left({\rm i}\omega^\ast\tau^\ast I-\int_{-1}^0d\eta_n(0,\theta){\rm e}^{{\rm i}\omega^\ast\tau^\ast\theta}\right)T_5^n=<\tilde F'''_{z z\bar z}+\frac{2{\rm i} g_{11}}{\omega^\ast\tau^\ast}\tilde F''_{z z}-\frac{{\rm i} \bar g_{02}}{3\omega^\ast\tau^\ast}\tilde F''_{\bar z\bar z}-\frac{{\rm i} (2\bar g_{11}+g_{20})}{\omega^\ast\tau^\ast}\tilde F''_{z\bar z},\beta_n>\\
~~~+\left[2g_{11}\frac{{\rm i} g_{20}}{\omega^\ast\tau^\ast}-(2\bar g_{11}+g_{20})\frac{{\rm i} g_{11}}{\omega^\ast\tau^\ast}-\frac{{\rm i} \bar g_{02} g_{02}}{3\omega^\ast\tau^\ast}+g_{21}\right]\left(\int_{-1}^0d\eta_n(0,\theta)\theta{\rm e}^{{\rm i}\omega^\ast\tau^\ast\theta}-I\right)<q(0)b_{n_0},\beta_n>,\end{array}$$
$$\begin{array}{l}\left(-{\rm i}\omega^\ast\tau^\ast I-\int_{-1}^0d\eta_n(0,\theta){\rm e}^{-{\rm i}\omega^\ast\tau^\ast\theta}\right)T_6^n=<\tilde F'''_{z \bar z\bar z}+\frac{{\rm i} g_{02}}{3\omega^\ast\tau^\ast}\tilde F''_{z z}-\frac{2{\rm i} \bar g_{11}}{\omega^\ast\tau^\ast}\tilde F''_{\bar z\bar z}+\frac{{\rm i} (2 g_{11}+\bar g_{20})}{\omega^\ast\tau^\ast}\tilde F''_{z\bar z},\beta_n>\\
~~~+\left[g_{02}\frac{{\rm i} \bar g_{02}}{3\omega^\ast\tau^\ast}+(2 g_{11}+\bar  g_{20})\frac{{\rm i} \bar  g_{11}}{\omega^\ast\tau^\ast}-\frac{2{\rm i} \bar g_{11}\bar g_{20}}{\omega^\ast\tau^\ast}+\bar g_{21}\right]\left(\int_{-1}^0d\eta_n(0,\theta)\theta{\rm e}^{-{\rm i}\omega^\ast\tau^\ast\theta}-I\right)<\bar q(0)b_{n_0},\beta_n>,\end{array}$$
$$\left(-3{\rm i}\omega^\ast\tau^\ast I-\int_{-1}^0d\eta_n(0,\theta){\rm e}^{-3{\rm i}\omega^\ast\tau^\ast\theta}\right)T_7^n=<\tilde F'''_{\bar z\bar z\bar z}+\frac{3{\rm i} \bar g_{20}}{\omega^\ast\tau^\ast}\tilde F''_{\bar z\bar z}+\frac{{\rm i} g_{02}}{\omega^\ast\tau^\ast}\tilde F''_{z\bar z},\beta_n>,$$
$$\begin{array}{l}\left(2{\rm i}\omega^\ast\tau^\ast I-\int_{-1}^0d\eta_n(0,\theta){\rm e}^{2{\rm i}\omega^\ast\tau^\ast\theta}\right)T_8^n=\left(3g_{21}+\frac{3{\rm i} g_{11} g_{20}}{\omega^\ast\tau^\ast}-\frac{6{\rm i} g_{11} \bar g_{11}}{\omega^\ast\tau^\ast}-\frac{{\rm i} g_{02} \bar g_{02}}{\omega^\ast\tau^\ast}\right)\\
\times\left(\int_{-1}^0d\eta_n(0,\theta){\rm e}^{2{\rm i}\omega^\ast\tau^\ast\theta}-I\right)T_1^n+<\tilde F'''_{z z z\bar z}+\frac{3{\rm i} g_{11}}{\omega^\ast\tau^\ast}\tilde F'''_{z z z}-\frac{3{\rm i} (g_{20}+\bar g_{11})}{\omega^\ast\tau^\ast}\tilde F''_{z z\bar z}-\frac{{\rm i} \bar g_{02}}{\omega^\ast\tau^\ast}\tilde F'''_{z\bar z\bar z\bar z},\beta_n>\\
+\left[\frac{6g_{11} \bar g_{11}}{(\omega^\ast\tau^\ast)^2}+\frac{15   g_{20} g_{11}}{(\omega^\ast\tau^\ast)^2}+\frac{ g_{02}\bar g_{02}}{3(\omega^\ast\tau^\ast)^2}\right]<\tilde F_{zz},\beta_n>
-\left[\frac{{\rm i} \bar g_{03} }{4\omega^\ast\tau^\ast}+\frac{3 \bar g_{02} g_{20}}{4(\omega^\ast\tau^\ast)^2}+\frac{ 5\bar g_{11}\bar g_{02}}{4(\omega^\ast\tau^\ast)^2}\right]<\tilde F_{\bar z\bar z},\beta_n>\\
+\left[\frac{7 \bar g_{02} g_{11}}{2(\omega^\ast\tau^\ast)^2}-\frac{9 \bar g_{11} g_{20}}{2(\omega^\ast\tau^\ast)^2}-\frac{\bar g_{02} \bar g_{20}}{2(\omega^\ast\tau^\ast)^2}-\frac{3 \bar g_{11}^2}{(\omega^\ast\tau^\ast)^2}-\frac{{\rm i}  g_{30}}{2\omega^\ast\tau^\ast}-\frac{3{\rm i} \bar g_{12}}{2\omega^\ast\tau^\ast} -\frac{3   g_{20}^2}{2(\omega^\ast\tau^\ast)^2}\right]<\tilde F_{z\bar z},\beta_n>\end{array}$$
and
$$\begin{array}{l}-\int_{-1}^0d\eta_n(0,\theta)T_9^n=2\left(\frac{2{\rm i} g_{11} g_{20}}{\omega^\ast\tau^\ast}-\frac{{\rm i} \bar g_{11} \bar g_{20}}{\omega^\ast\tau^\ast}+g_{21}+\bar g_{21}\right)
\left(\int_{-1}^0d\eta_n(0,\theta)\theta-I\right)T_2^n\\
+<\tilde F''''_{z z \bar z\bar z}+\frac{{\rm i} g_{02}}{3\omega^\ast\tau^\ast}\tilde F'''_{z z z}-\frac{{\rm i}\bar g_{02}}{3\omega^\ast\tau^\ast}\tilde F'''_{\bar z \bar z\bar z}+\frac{{\rm i} (\bar g_{20}+4g_{11}}{(\omega^\ast\tau^\ast)^2}\tilde F'''_{z z\bar z}-\frac{{\rm i} ( g_{20}+4\bar g_{11}}{(\omega^\ast\tau^\ast)^2}\tilde F'''_{z \bar z\bar z},\beta_n>\\
+\left[\frac{{\rm i} g_{12}}{\omega^\ast\tau^\ast}-\frac{4   g_{11}^2}{(\omega^\ast\tau^\ast)^2}-\frac{ g_{11}\bar g_{20}}{(\omega^\ast\tau^\ast)^2}+\frac{ 2\bar g_{11}  g_{02}}{(\omega^\ast\tau^\ast)^2}\right]<\tilde F_{zz},\beta_n>\\
+\left[-\frac{{\rm i} \bar g_{12} }{\omega^\ast\tau^\ast}-\frac{4 \bar g_{11}^2}{(\omega^\ast\tau^\ast)^2}-\frac{ \bar g_{11}  g_{20}}{(\omega^\ast\tau^\ast)^2}+\frac{2 \bar g_{02}  g_{11}}{(\omega^\ast\tau^\ast)^2}\right]<\tilde F_{\bar z\bar z},\beta_n>\\
+2\left[\frac{3   g_{20} g_{11}}{(\omega^\ast\tau^\ast)^2}+\frac{3 \bar g_{11} \bar g_{20}}{(\omega^\ast\tau^\ast)^2}+\frac{8 |g_{11}|^2}{(\omega^\ast\tau^\ast)^2}+\frac{ |g_{20}|^2}{(\omega^\ast\tau^\ast)^2}+\frac{| g_{02}|^2}{3(\omega^\ast\tau^\ast)^2}\right]<\tilde F_{z\bar z},\beta_n>.\end{array}$$
Recalling Eq.(\ref{FFF}),  all the variables $T_j$'s can be solved explicitly, so we can determine the first and second Lyapunov coefficients at any $(k,\tau^\ast)$ by using Eqs.(\ref{l1l1}-\ref{l2l2}).





\end{document}